\documentclass[12pt]{article}
\usepackage{amscd,amsfonts,amssymb,amsmath,amsthm,latexsym,mathrsfs,graphicx,setspace}
\usepackage[OT2,T1]{fontenc}
\DeclareSymbolFont{cyrletters}{OT2}{wncyr}{m}{n}

\newcommand{\alg}{\text{alg}}

\newcommand{\bbF}{\mathbb{F}}

\newcommand{\bbGa}{[\![\Ga]\!]}

\newcommand{\bbQ}{\mathbb{Q}}

\newcommand{\bbZ}{\mathbb{Z}}

\newcommand{\bbZp}{{\mathbb{Z}_p}}

\newcommand{\bfD}{\mathbf{D}}

\newcommand{\bfL}{\mathbf{L}}

\newcommand{\bfR}{\mathbf{R}}
\newcommand{\bfT}{\mathbf{T}}

\newcommand{\bra}[1]{{\langle #1 \rangle}}

\newcommand{\calI}{\mathcal{I}}

\newcommand{\calL}{\mathcal{L}}

\newcommand{\calO}{\mathcal{O}}
\newcommand{\calS}{\mathcal{S}}

\DeclareMathOperator{\chr}{char}
\DeclareMathOperator{\coker}{coker}
\DeclareMathOperator{\Cone}{Cone}
\DeclareMathOperator{\corank}{corank}
\newcommand{\cn}{\colon}

\newcommand{\cycl}{\text{cycl}}

\renewcommand{\emptyset}{\varnothing}
\newcommand{\ep}{\epsilon}

\newcommand{\fkh}{\mathfrak{h}}

\newcommand{\fkm}{\mathfrak{m}}

\newcommand{\fkp}{\mathfrak{p}}

\DeclareMathOperator{\Frac}{Frac}
\newcommand{\Frob}{\text{Frob}}
\newcommand{\Ga}{\Gamma}
\newcommand{\ga}{\gamma}
\DeclareMathOperator{\Gal}{Gal}
\DeclareMathOperator{\GL}{GL}

\DeclareMathOperator{\Hom}{Hom}

\DeclareMathOperator{\img}{img}
\newcommand{\inv}{^{-1}}

\newcommand{\La}{\Lambda}
\newcommand{\la}{\lambda}
\newcommand{\llim}{\varprojlim}

\newcommand{\Nekovar}{{Nekov\'a\v r}}

\DeclareMathOperator{\ord}{ord}

\newcommand{\ov}[1]{{\overline{#1}}}

\DeclareMathOperator{\rank}{rank}

\newcommand{\rlim}{\varinjlim}

\newcommand{\scrD}{\mathscr{D}}

\DeclareMathSymbol{\Sha}{\mathalpha}{cyrletters}{"58}
\DeclareMathOperator{\SL}{SL}
\DeclareMathOperator{\Spec}{Spec}

\newcommand{\sto}{^*(1)}

\newcommand{\wh}[1]{{\widehat{#1}}}
\newcommand{\wt}[1]{{\widetilde{#1}}}

\newcommand{\ddcoft}{\text{-}\mathbf{coft}}
\newcommand{\ddfl}{\text{-}\mathbf{fl}}
\newcommand{\ddft}{\text{-}\mathbf{ft}}
\newcommand{\ddmod}{\text{-}\mathbf{mod}}

\newtheorem{thm}{Theorem}[section]
\newtheorem{ppn}[thm]{Proposition}
\newtheorem{lem}[thm]{Lemma}
\theoremstyle{definition}

\theoremstyle{remark}
\newtheorem{rem}[thm]{Remark}

\author{Jonathan Pottharst}
\title{Parity-induced Selmer Growth for Symplectic, Ordinary Families}

\begin{document}

\maketitle

\abstract{Let $p$ be an odd prime, and let $K/K_0$ be a quadratic
  extension of number fields.  Denote by $K_\pm$ the maximal
  $\bbZ_p$-power extensions of $K$ that are Galois over $K_0$, with
  $K_+$ abelian over $K_0$ and $K_-$ dihedral over $K_0$.  In this
  paper we show that for a Galois representation over $K_0$ satisfying
  certain hypotheses, if it has odd Selmer rank over $K$ then for one
  of $K_\pm$ its Selmer rank over $L$ is bounded below by $[L:K]$ for
  $L$ ranging over the finite subextensions of $K$ in $K_\pm$.  Our
  method of proof generalizes a method of Mazur--Rubin, building upon
  results of \Nekovar, and applies to abelian varieties of arbitrary
  dimension, (self-dual twists of) modular forms of even weight, and
  (twisted) Hida families.}


\tableofcontents

\section{Introduction}\label{sect-intro}

This paper fits into a circle of ideas that might broadly be called
arithmetic in $p$-adic dihedral extensions of number fields.  A
stunning result in this area follows from the formula of Gross--Zagier
\cite{gz} and Kolyvagin's Euler system \cite{koly}.  These works
establish the existence of a wealth of points on elliptic curves in
certain dihedral extensions of $\bbQ$.  More precisely, let $E/\bbQ$
be a (modular) elliptic curve of conductor $N$ without complex
multiplication, and let $K/\bbQ$ be an imaginary quadratic field
satisfying the Heegner condition relative to $N$.  For a prime $p>3$
that is split in $K$ and at which $E$ has good, ordinary reduction, we
let $K_-/K$ denote the $p$-anticyclotomic extension.  Suppose that
$E(K)$ has rank one; then for every finite subextension $L$ of
$K_-/K$, one has $\rank_\bbZ E(L) \geq [L:K]$.  (The Heegner points
supply the desired rank.)

We call $K_-/\bbQ$ a {\it $p$-adic dihedral extension}, since it is
Galois; it contains a quadratic subextension $K/\bbQ$ with $K_-/K$ a
$\bbZ_p$-power extension; and the conjugation action of (any lift to
$K_-$ of) the nontrivial automorphism of $K/\bbQ$ upon $\Gal(K_-/K)$
is by means of inversion.  Given a self-dual motive over $\bbQ$, a
general yoga of signs of functional equations leads us to expect that
if the underlying geometric object has an odd number of rational
cycles over $K$, then it acquires rational algebraic cycles in number
at least equal to $[L:K]$, for $L$ ranging over the finite
subextensions of $K_-/K$.  In the case considered by Gross, Zagier,
and Kolyvagin, this can be seen on the levels of $L$-functions,
algebraic cycles, and Selmer groups all at once.

Our present work generalizes a weakened form of the above result.  Let
$p$ be an odd prime, and let $K/K_0$ be a quadratic extension of
number fields.  Denote by $K_\pm$ the maximal $\bbZ_p$-power
extensions of $K$ that are Galois over $K_0$, with $K_+$ abelian over
$K_0$ and $K_-$ $p$-adic dihedral over $K_0$.  We show that for a
Galois representation over $K_0$ satisfying certain hypotheses, if it
has odd Selmer rank over $K$ then for one of $K_\pm$ its Selmer rank
over $L$ is bounded below by $[L:K]$ for $L$ ranging over the finite
subextensions of $K$ in $K_\pm$.  (See \S\ref{sub-cohom} for the
precise statement.)

Our method of proof generalizes a method of Mazur--Rubin, building
upon results of \Nekovar, and applies to abelian varieties of
arbitrary dimension, (self-dual twists of) modular forms of even
weight, and (twisted) Hida families.  It should be noted that
\Nekovar\ has obtained a much more precise result (cf.\
\cite[10.7.15.iii]{nekovar}) by making use of the Cassels--Tate
pairing apparatus, which we do not need.

The work falls naturally into two parts.  The first part is a purely
cohomological result for $p$-adic Galois representations; we just say
here that it requires a self-dual, ``ordinary'' representation with no
Tamagawa obstructions, no exceptional zeroes, and a weak form of
residual irreducibility.  The representation may be free over any
complete, Noetherian, Gorenstein local ring with finite residue field.
The second part of our work is to show when the theorem applies to
concrete examples.

The main insight is that by working only in the cohomological arena,
we need not appeal to any deep conjectures to prove our result; this
allows us to obtain very general evidence for the framework of
conjectures.  On the other hand, {\it any} application of the result
to producing algebraic cycles requires the resolution of these
conjectures.  We also underscore the obvious drawback of our approach:
it is unable to distinguish between $K_+$ and $K_-$!

There are quite diverse circumstances under which the ability to
descend a variety below a given base field conjecturally implies the
existence of rational cycles.  For example, in \cite{djk}, it is shown
how to produce pieces in the $\ell$-adic \'etale cohomology of certain
varieties over a finite field which, via Tate's conjectures, predict
the existence of rational algebraic cycles ``for free'' after
enlarging the finite field.

This paper has benefited from our interactions with numerous people.
We especially thank Barry Mazur for suggesting we take up this
project, and for his advice and encouragement throughout; and Jan
\Nekovar\ for his arranging our study in Paris and Luminy, for his
thorough reading of an earlier draft and helpful conversations
thereafter, and not least for his investing the time to write out
\cite{nekovar}.  Matthew Emerton brought \cite{epw} to our attention,
which sprang forth to help at a critical juncture.  Finally, we
heartily thank the NSF for its support through the MSPRFP, under which
this work was produced, and l'Institut de Math\'ematiques de Jussieu
for its hospitality.

To conclude the introduction, we describe the contents of the paper.
Immediately following, we begin with our notation and running
hypotheses.  In the second section, we prove our main technical
result, using the methods of Mazur--Rubin and Nekov\'a\v r.  In the
third section, we deduce applications to the Selmer groups of abelian
varieties, modular forms, and Hida families.  In the fourth section,
we give concrete numerical examples.  Finally, in an appendix, we give
a brief exposition to those aspects of Nekov\'a\v r's Selmer complexes
that will be of use to us.

\subsection{Notation and hypotheses}\label{sub-notation}

The following notation will be in force through this paper.

$p$ denotes an {\it odd} rational prime.

$(R,\fkm,k)$ is a Gorenstein, complete, Noetherian local ring, with
finite residue field of characteristic $p$, and $d = \dim R$.  For any
abelian profinite group $G$, the ring $R[\![G]\!]$ is the completed
group algebra of $G$ with coefficients in $R$.  For $g \in G$, we
write $\bra{g}$ for the corresponding grouplike element of
$R[\![G]\!]$.

If $K/K_0$ is a quadratic extension of finite extensions of $\bbQ$, we
let $\ov{K}$ denote a fixed algebraic closure of $K$, with Galois
group $G_K = \Gal(\ov{K}/K)$.  Denote the maximal
$\bbZ_p$-power-extension of $K$ in $\ov{K}$ by $K_\infty$.  For any
subextension $L$ of $K_\infty/K$ whatsoever, we set $\Ga_L :=
\Gal(L/K)$ and $\La_L := R[\![\Ga_L]\!]$, and write $\calI_L :=
\ker(\La_{K_\infty} \to \La_L)$ for its relative augmentation ideal.
Moreover, for any $\La_{K_\infty}$-module $M$, we set $M_L := M
\otimes_\La \La_L = M/\calI_L M$.  We drop the subscripts when $L =
K_\infty$ (since $\calI_{K_\infty}=0$ anyway).  The involution $\iota$
given by inversion $\ga \mapsto \ga\inv$ on $\Ga$ induces an action on
all groups and group rings above, and we denote these actions all by
the same symbol.  For a $\La$-module $M$, if we twist the $\La$-module
structure through $\iota$, the result is written $M^\iota$.

Denote by $\sigma \in \Gal(K/K_0)$ the nontrivial element.  Any choice
of lift of $\sigma$ to $\Gal(K_\infty/K_0)$ acts by conjugation on
$\Ga$.  This allows us to decompose $\Ga$ as a product $\Ga_+ \times
\Ga_-$ of eigenspaces, so that (any lift of) $\sigma$ acts on
$\Ga_\pm$ via $\ga \mapsto \ga^{\pm1}$.  We let $K_\pm$ be the fixed
field of $\Ga_\mp$, so that $\Ga_\pm$ is naturally identified with
$\Ga_{K_\pm}$.  In general, simply write ``$\pm$'' in subscripts in
place of ``$K_\pm$'', as in $\La_\pm$, $\calI_\pm$, and $M_\pm$.

For each place $v$ of $K$ we fix an algebraic closure $\ov{K}_v$ of
$K_v$, with Galois group $G_v = \Gal(\ov{K}_v/K_v)$.  We also fix
embeddings $\ov{K} \hookrightarrow \ov{K}_v$, which induce inclusions
$G_v \hookrightarrow G_K$ as decomposition groups.  If $I_v \subset
G_v$ is the inertia subgroup, we write $\Frob_v \in G_v/I_v$ for the
arithmetic Frobenius element.

By $S$ we mean a finite set of places of $K_0$, assumed to contain all
Archimedean places and all places lying over $p$.  The set $S_f$
consists of all finite places of $S$, and we partition $S_f = \Sigma
\sqcup \Sigma'$, with $v \in \Sigma$ if and only if $v$ lies over $p$.
We let $K_S$ be the maximal extension of $K$ in $\ov{K}$ unramified
outside $S$, with Galois group $G_{K,S} = \Gal(K_S/K)$.  (Note that
$K_\infty \subseteq K_S$ by \cite[\S5.4, Lemma (i)]{lang}.)  For each
place $v$ of $K$, our choices provide us with a composite map $G_v
\hookrightarrow G_K \twoheadrightarrow G_{K,S}$.

For any ring $A$, we write $\bfD(A)$ for the derived category of
$A$-modules.  By $\bfD_\text{perf}^{[a,b]}(A)$, where $a \leq b$ are
integers, we mean the subcategory of $\bfD(A)$ whose objects can be
represented by complexes $C^\bullet$ with the $C^i$ finitely generated
and free over $A$, and nonzero only for $a \leq i \leq b$.

We write $D$ for the Pontryagin duality functor
$\Hom_{\bbZ_p}(\cdot,\bbQ_p/\bbZ_p)$, and $\scrD$ for the functor
$\Hom_R(\cdot,R)$.  The latter operation will {\it only} be applied to
complexes $T^\bullet$ consisting of {\it free} $R$-modules; and in
this scenario, we consider $\scrD(T^\bullet)$ as providing an explicit
choice of complex representing the Grothendieck dual of $T^\bullet$ in
$\bfD(R)$.  (See \S\ref{sub-inter} for more on these matters.)  We
also let $\scrD_\La$ mean $\Hom_\La(\cdot,\La)$, with the same
proviso.

We denote by $\chi_\cycl \cn G_\bbQ \to \bbZ_p^\times$ the $p$-adic
cyclotomic character.  As $\bbZ_p^\times$ decomposes canonically as
$\mu_{p-1} \times (1+p\bbZ_p)$, we accordingly write $\chi_\cycl =
\tau_\cycl \cdot \ga_\cycl$.  The Galois character $\tau_\cycl$
corresponds via class field theory to the unique Dirichlet character
of conductor $p$ whose reduction modulo $p$ is the mod $p$ cyclotomic
character; we abusively use $\tau_\cycl$ to denote this Dirichlet
character as well.  For later use, we note that, since $1+p\bbZ_p$ is
free of rank one over $\bbZ_p$, the character $\ga_\cycl$ admits a
unique square root.  On the other hand, since $\mu_{p-1}$ is cyclic,
the character $\tau_\cycl^i$ only admits a square root when $i$ is
even, in which case there are two, corresponding to the lifts of $i
\in \bbZ/(p-1)$ to $\bbZ/2(p-1)$.

If $\chi \cn G_K \to R^\times$ is a character, we write $R(\chi)$ for
a free rank one $R$-module on which $G_K$ acts through multiplication
by $\chi$.  And if $M$ is an $R[G_K]$-module we write $M(\chi)$ for $M
\otimes_R R(\chi)$, with $G_K$ acting diagonally.  As a shorthand, we
let $M(n)$ denote $M(\chi_\cycl^n)$ for $n \in \bbZ$.

Given an $R$-module $M$, one can associate to each minimal prime $\fkp
\in \Spec(R)$ the {\it rank} of $M$ at $\fkp$.  Namely, the ring
$T = (R_\fkp)^\text{red}$ is a domain, hence
\[
\rank_{R,\fkp} M := \rank_T (M \otimes_R T) := \dim_{\Frac T} (M
\otimes_R \Frac T)
\]
makes sense.  We think of $\rank_R M$ as the system $(\rank_{R,\fkp}
M)_\fkp$ of nonnegative integers, with $\fkp$ ranging over the minimal
primes in $\Spec(R)$.  (One similarly defines $\corank = \rank \circ
D$.)  In particular, if $M \approx R^{\oplus d}$ is free of rank $d$
(in the usual parlance), then $\rank_R M$ has constant value $d$.
Moreover, to say ``$\rank_R M$ is even (resp.\ odd, $\geq N$)'' means
that the claim holds for each $\rank_{R,\fkp} M$.  The ring $\La$
satisfies the same hypotheses as $R$ does, and all the above applies
to ``$\rank_\La M$'' as well.  On the other hand, the natural
homomorphisms $R \to \La \to \La/\calI=R$ induce maps $\Spec(R) \to
\Spec(\La) \to \Spec(R)$ that are mutually inverse {\it bijections} on
minimal primes (apply \cite[8.9.7.i]{nekovar} with $\Delta=\{1\}$),
which we treat as an identification.  For $M$ a finite free
$\La$-module, one clearly has $\rank_\La M = \rank_R M \otimes_\La R$
under this identification.

\section{Technical results}

In this section we prove our cohomological theorem on Selmer growth.

The larger body of results we draw upon comes from \Nekovar's
formalism of Selmer complexes, which expresses the arithmetic local
and global dualities in the language of derived categories.  Since
familiarity with these ideas is not necessary for our proofs, we cite
the relevant theorems from \cite{nekovar}.  The basic constructions
are sketched in an appendix, for the benefit of the reader who wishes
to know of their origins.

We also use two key ideas of Mazur--Rubin that allow one to force
algebraic $p$-adic $L$-functions to have zeroes, granted they obey
certain functional equations.  In the earlier portion of this section,
we begin by recalling Mazur--Rubin's ideas in the appropriate
generality.  In the latter part, we prove our cohomological theorem.

\subsection{Skew-Hermitian complexes and functional equations}
\label{sub-mr}

In this section we review the method of Mazur--Rubin, stating their
results in a generality that is suitable for our needs.

Consider a complex $C^\bullet = [\Phi \stackrel{u}{\hookrightarrow}
\Psi]$, concentrated in degrees $[1,2]$, with $\Phi,\Psi$ finite free
over $\La$ of the same rank, and $u$ injective.  Assume that this
complex is equipped with a quasi-isomorphism $\alpha \cn C^\bullet
\stackrel{\sim}{\to} \Hom(C^\bullet,\La)^\iota[-3]$ satisfying
$\Hom(\alpha,\La)^\iota[-3] = -\alpha$ up to chain homotopy.  The
following provides an example of such a complex.

Let $M$ be free of finite rank over $\La$, equipped with a
nondegenerate, skew-Hermitian $\La$-bilinear pairing $h\cn M \otimes
M^\iota \to \La$, with image contained in $\fkm$.  If we write $M^* :=
\Hom_\La(M^\iota,\La)$, the adjoint $h^\text{ad} \cn M \to M^*$ serves
as the boundary operator of a complex $[M \stackrel{h^\text{ad}}{\to}
M^*]$ concentrated in degrees $1,2$; the nondegeneracy of $h$ means
that $h^\text{ad}$ is injective.  This complex is equipped with an
obvious duality pairing.  The complex just described, together with
its duality structure, is denoted $C(M,h)^\bullet$ and called a {\it
basic skew-Hermitian complex}.

We will make use of the following two propositions of Mazur--Rubin.

\begin{ppn}\label{mr-rep}
Every $C^\bullet$ is quasi-isomorphic to a $C(M,h)^\bullet$, in a
manner respecting the duality pairings.
\end{ppn}
\begin{proof}
This is \cite[Proposition 6.5]{mr:org}.
\end{proof}

The proof of existence relies crucially on Nakayama's lemma, and thus
on the fact that $R$ is local.  The author sees no means to generalize
the methods of \cite{mr:org} beyond the local case.

\begin{ppn}\label{mr-gen}
Let $\Xi$ be a finite group of commuting involutions of $\Ga$.
Suppose $I \subset \La$ is a nonzero principal ideal that is preserved
by $\Xi$.  Then there exists a generator $\calL \in I$ such that
$\calL^\xi = \ep(\xi)\calL$ for some homomorphism $\ep \cn \Xi \to
\{\pm1\}$.  The value $\ep(\xi)$ depends only on $\xi$ and $I$, and
not on $\Xi$ or $\calL$.
\end{ppn}
\begin{proof}
The proof in \cite[Proposition 7.2]{mr:growth} refers to the case
where $R$ is the ring of integers in a finite extension of $\bbQ_p$,
but it applies without change to any $R$ under our hypotheses.
\end{proof}

Consider a basic skew-Hermitian complex $C^\bullet = C(M,h)^\bullet$
over $\La$, as above.  We will apply the preceding proposition to the
characteristic ideal
\[
I := \chr_\La H^2(C^\bullet) = \chr_\La(\coker h^\text{ad}) =
\det(h^\text{ad})\La
\]
(if it is nonzero), and to the group $\Xi$ generated by the two
involutions $\iota$ and $\sigma$, which we now recall.

First, one always has the inversion involution $\iota \cn \ga \mapsto
\ga\inv$.  We can use the skew-Hermitian property of $h$ to calculate
that
\[
\det(h^\text{ad})^\iota = \det(h^{\text{ad}\,\iota}) =
\det(-h^\text{ad}) = (-1)^r\det(h^\text{ad}),
\]
with $r := \rank_\La M \pmod{2}$.  This shows that $I$ is stable under
$\iota$, and moreover that $\ep(\iota) = (-1)^r$.  We point out that
$r$ may be computed ``over $R$'' as follows.  The involution $\iota$
acts trivially on $R$, so $h^\text{ad} \pmod{\calI}$ is
skew-symmetric; therefore, $\rank_{R,\fkp} \left(\img h^\text{ad}
\otimes_\La R\right)$ is even for every minimal prime $\fkp$ of $R$.
This lets us calculate that, for all such $\fkp$,
\[\begin{split}
r &= \rank_\La M^*
= \rank_{R,\fkp} M^* \otimes_\La R \\
&\equiv \rank_{R,\fkp}
  (M^* \otimes_\La R) / (\img h^\text{ad} \otimes_\La R) \pmod{2} \\
&= \rank_{R,\fkp}
  (\coker h^\text{ad} \otimes_\La R),
\end{split}\]
where the last equality is by the right exactness of $\otimes$.

For the other involution, we recall that we are given a degree $2$
subfield $K_0$ of $K$ as in \S\ref{sub-notation}, and we have the
involution $\sigma$ which acts on $\Ga_\pm$ via $\ga \mapsto
\ga^{\pm1}$.  In the next section, our complex $C$ will arise
functorially from a Galois module $T$ defined over $K_0$.  Each lift
of $\sigma$ to $\Gal(K_\infty/K_0)$ will induce an isomorphism $C
\stackrel{\sim}{\to} C^\sigma$, and therefore $I = \chr_\La H^2(C)$
will be stable under $\sigma$.

\subsection{The cohomological theorem}\label{sub-cohom}

We begin our discussion of our main theorem on Selmer growth by laying
out the setup and hypotheses.

Continue with notations as in \S\ref{sub-notation}.  Let $T$ be a
nonzero, free, finite rank $R$-module with a continuous, linear
$G_{K_0,S}$-action.  We require the following list of hypotheses and
data attached to $T$ (whose motivations are explained in Remark
\ref{rem-hyp-explain}):

\begin{description}

\item[(Symp)] $T$ is symplectic; i.e., it is equipped with a
Galois-equivariant perfect pairing $T \otimes T \to R(1)$, and hence
an isomorphism $j \cn T \stackrel{\sim}{\to} \scrD(T)(1)$, that is
skew-symmetric in the sense that $\scrD(j)(1) = -j$.

\item[(Ord)] For each $v \in \Sigma$, we are given a $G_v$-stable
$R$-direct summand $T_v^+ \subset T$ that is Lagrangian for the
symplectic structure: $j(T_v^+) = \scrD(T/T_v^+)(1)$.  Set $T_v^- =
T/T_v^+$, obtaining an exact sequence:
\[
0 \to T_v^+ \to T \to T_v^- \to 0.
\]

\item[(Tam)] For every place $v$ of $K$ lying over $\Sigma'$, the
submodule $T^{I_v}$ is assumed {\it free} over $R$.  Moreover, for
such $v$, the operator $\Frob_v-1$ acts bijectively on $H^1(I_v,T)$.

\end{description}

Set $A = D(\scrD(T)) \cong D(T)(1)$ (by the the self-duality of $T$)
and $A_v^\pm = D(T_v^\mp)(1) \subset A$, obtaining an exact sequence:
\[
0 \to A_v^+ \to A \to A_v^- \to 0.
\]
One can view $A$ as the ``divisible incarnation'' of $T$.  Specifying
$T_v^+$ is equivalent to specifying $A_v^+$.  For an extension $L$ of
$K$ and a place $v$ of $L$ lying over $v_0 \in \Sigma$, we also set
$T_v^{\pm} = T_{v_0}^{\pm}$ and $A_v^\pm = A_{v_0}^\pm$.

Here are our two remaining hypotheses:

\begin{description}

\item[(Irr)] The following morphism (in which $A^{G_{K,S}}$ maps
diagonally) is injective.
\[
A^{G_{K,S}} \to \bigoplus_{v \text{ lying over } \Sigma} A_v^-
\]

\item[(Zero)] For all places $v$ of $K$ lying over $\Sigma$, one has
$(T_v^-/\fkm)^{G_v} = 0$.

\end{description}

Given the above data, we may define the {\it (strict Greenberg) Selmer
groups} of $A$ over any algebraic extension $L$ of $K$.  For each
prime $v$ of $L$, choose a decomposition group $G_v$ above $v$ in
$\Gal(\ov{K}/L)$ and let $I_v$ be its inertia subgroup.  Then the
Selmer group of $A$ over $L$ is
\begin{equation}\label{eqn-sel-defn}
S_A(L) := \ker\left[ H^1(\Gal(\ov{K}/L),A) \to \prod_{v \nmid p}
  H^1(I_v,A) \times \prod_{v \mid p} H^1(G_v,A_v^-) \right].
\end{equation}

We now state our main technical theorem.

\begin{thm}\label{thm-cohom}
With notation as in \S\ref{sub-notation}, assume the above five
hypotheses hold for $T$.  If $S_A(K)$ has odd $R$-corank, then for at
least one choice of sign $\ep = \pm$, we have $\corank_R S_A(L) \geq
[L:K]$ for every finite extension $L/K$ contained in $K_\ep$.
\end{thm}

\begin{rem}\label{rem-hyp-explain}
The assumption (Ord) is a variant of assuming that $T$ is ``ordinary''
(or better: ``Pan\v{c}i\v{s}kin'') above $p$, and (Tam) is related to
$p$ not dividing a ``Tamagawa number'' at $v$, for all $v \nmid p$;
see Remark \ref{rem-DVR-simplify} below.  The condition (Irr) holds,
in particular, if $A[\fkm]^{G_{K,S}} = 0$, and {\it a fortiori} if
$A[\fkm]$ is an irreducible residual representation (note that its
rank is at least two).  The hypothesis (Zero) excludes the case of an
``algebraic exceptional zero'' playing a similar role to those found
by Mazur--Tate--Teitelbaum in \cite{mtt}.
\end{rem}

\begin{rem}\label{rem-DVR-simplify}
In the case where $R$ is the ring of integers $\calO$ in a finite
extension $F$ of $\bbQ_p$, some simplifications are possible.  First,
$\calO$ is a DVR, so that $A$ is isomorphic to $T \otimes_\calO
F/\calO$.  Second, $\calO$ is a PID, so that in (Ord), to determine
$T_v^+$, it suffices to specify the subspace $T_v^+ \otimes_\calO F
\subset V$, where $V = T \otimes_\calO F$; moreover, in (Tam), the
freeness assumption is automatically satisfied.  Third, by the
equation
\[
\scrD(\text{Err}_v^{ur}(\scrD,T)) \stackrel{\sim}{\longrightarrow}
D(\text{Err}_v^{ur}(\Phi,T))
\]
appearing in the proof of \cite[7.6.7.ii]{nekovar}, combined with the
calculation of \cite[7.6.9]{nekovar}, one can rephrase (Tam) as the
claim that for every $v \in \Sigma'$ one has
$H^1(I_v,T)_\text{tors}^{\Frob_v=1} = 0$.  The order of
$H^1(I_v,T)_\text{tors}^{\Frob_v=1}$ is precisely the ($p$-part of
the) {\it Tamagawa number} of $T$ at $v$ (see \cite[Proposition
I.4.2.2.ii]{fpr}).  The same computation also shows that in order to
verify this last criterion, it would suffice to show that
$V^{I_v}/T^{I_v} \twoheadrightarrow A^{I_v}$.
\end{rem}

\begin{rem}
By \cite[10.7.15.iii]{nekovar}, under almost identical hypotheses, one
can take $\ep = -1$, and for $L/K$ contained in $K_-$ one knows that
$\corank_R S_A(L)$ is odd.
\end{rem}

Our proof of the theorem relies on the following results of \Nekovar,
which allow us to represent $S_A(L)$ in terms of a certain Selmer
complex (cf.\ Proposition \ref{ppn-control}).  For all our Selmer
complexes, over any number field containing $K_0$, we choose the local
conditions $\Delta$ to be ``unramified'' at primes $v$ lying over
$\Sigma'$, and ``(strict) Greenberg'' at primes $v$ lying over
$\Sigma$.

\begin{lem}\label{lem-dual-cond}
The $\scrD(1)$-dual local conditions to $\Delta$ are isomorphic to
$\Delta$ under the identification $T \cong \scrD(T)(1)$.
\end{lem}
\begin{proof}
Plug (Tam) into \cite[7.6.12]{nekovar}, and plug (Symp) and (Ord) into
\cite[6.7.6.iv]{nekovar}.
\end{proof}

This proposition is \Nekovar's Iwasawa-theoretic arithmetic duality
(and perfectness) theorem:

\begin{ppn}\label{ppn-bound}
The Iwasawa-theoretic Selmer complex
\[
C := \wt{\bfR\Ga}_{f,\mathrm{Iw}}(K_\infty/K,T;\Delta)
\]
defined in \cite[8.8.5]{nekovar} lies in
$\bfD_\mathrm{perf}^{[1,2]}(\La)$, and is equipped with a
skew-Hermitian duality quasi-isomorphism
\[
\alpha \cn C \stackrel{\sim}{\longrightarrow} \scrD_\La(C)^\iota[-3].
\]
When we represent $C$ by a complex of the form $[\Phi
\stackrel{u}{\to} \Psi]$ with $\Phi,\Psi$ finite free over $\La$ (with
respective degrees $1,2$), the modules $\Phi,\Psi$ have the same rank.
Moreover, the $\La$-module
\[
\calS := H^2(C)
\]
is $\La$-torsion if and only if the differential $u$ is injective.
\end{ppn}
\begin{proof}
All the claims follow by copying the steps of \cite[9.7]{nekovar}
word-for-word.  In particular: To see that $C$ lies in
$\bfD_\text{perf}^{[0,3]}$, use the proof of \cite[9.7.2.ii]{nekovar},
noting the necessity of the freeness condition in (Tam) (which is
automatic in the case under \Nekovar's consideration).  That $\alpha$
is an isomorphism follows from Lemma \ref{lem-dual-cond} and
\cite[9.7.3.iv]{nekovar}.  It is skew-Hermitian by
\cite[9.7.7(ii)]{nekovar}.  The placement $C \in
\bfD_\text{perf}^{[1,2]}(\La)$ follows from (Irr) and
\cite[9.7.5.ii]{nekovar}.  The final two claims are
\cite[9.7.7.iv]{nekovar}.
\end{proof}

The letter ``$\calS$'' is meant to remind us of the word ``Selmer''.
This choice of mnemonic is because of the following comparison.

\begin{ppn}\label{ppn-control}
Let $S_A(K_\infty)$ be as in Equation \ref{eqn-sel-defn}.  For each
(possibly infinite) subextension $L$ of $K_\infty/K$, recalling that
$\calS_L := \calS \otimes_\La \La_L$, we have
\[
\calS_L \cong D(S_A(L))^\iota,
\qquad \text{i.e.} \qquad
S_A(L) \cong D(\calS_L)^\iota.
\]
\end{ppn}
\begin{proof}
We use the notation of \cite[9.6]{nekovar}, with the exception that
our $S_A(L)$ is written $S_A^\text{str}(L)$ there.  By
\cite[9.6.3]{nekovar}, for any $L$ as in the proposition, there is a
surjection $\wt{H}_f^1(L/K,A) \twoheadrightarrow S_A(L)$, which is an
isomorphism provided that for all places $v \in \Sigma$, we have
$(A_v^-)^{G_v \cap G_L} = 0$.  By \cite[9.6.6.iii]{nekovar}, it
suffices to check the latter condition when $L=K$; by Nakayama's
lemma, this is equivalent to requiring that $(T_v^-/\fkm)^{G_v} \cong
(A_v^-[\fkm])^{G_v} = 0$, which is precisely (Zero).

Let $L$ be any subextension of $K_\infty/K$.  Invoking
\cite[9.7.2.i]{nekovar}, we find that
\begin{equation}\label{eqn-control}
D(S_A(L))^\iota \cong H^2(C(L)),
\end{equation}
where $C(L)$ is the Iwasawa-theoretic Selmer complex constructed just
like $C$ was, but with $L$ in place of $K_\infty$.  In particular, our
proposition follows in the case $L = K_\infty$.

We now invoke \Nekovar's control theorem \cite[8.10.10]{nekovar} (cf.\
the discussion at the end of \S\ref{sub-global}), showing that $C(L)
\cong C \otimes^\bfL_\La \La_L$.  Represent $C$ by a complex
$C^\bullet$ of the form $[\Phi \stackrel{u}{\to} \Psi]$ as in
\ref{ppn-bound}.  Since $\Phi,\Psi$ are free, the object $C(L)$ is
represented by the complex $C^\bullet \otimes_\La \La_L$.  Therefore,
\[
H^2(C(L)) \cong H^2(C^\bullet \otimes_\La \La_L) = \coker(u \bmod
\calI_L) = \coker(u) \bmod \calI_L = \calS_L,
\]
which, together with Equation \ref{eqn-control}, proves the
proposition in general.
\end{proof}

In particular, under our hypotheses, a form of ``perfect control''
holds: the natural maps $S_A(L) \to S_A(L')^{\Gal(L'/L)}$ are
isomorphisms, for any $K_\infty/L'/L/K$.

\begin{proof}[Proof of Theorem \ref{thm-cohom}]
Recall that when $L = K_\pm$, we write ``$\pm$'' as a subscript
instead of ``$K_\pm$''.  It suffices to show that at least one of
$\calS_\pm$ is not torsion over its respective $\La_\pm$, because then
for every finite subextension $L/K$ of $K_\pm$, we have
\[
\corank_R S_A(L)
= \rank_R (\calS_\pm \otimes_{\La_\pm} \La_L)
\geq \rank_{\La_\pm} \calS_\pm \cdot \rank_R \La_L \geq 1 \cdot [L:K],
\]
as was desired.  If $\calS$ is not a torsion $\La$-module, then both
of the $\calS_\pm$ are not torsion $\La_\pm$-modules, and our theorem
follows trivially.  So let us assume henceforth that $\calS$ is
torsion over $\La$.  In this case, the characteristic ideal $\chr_\La
\calS \subseteq \La$ is nonzero, and $(\chr_\La \calS)\La_\pm$ divides
$\chr_{\La_\pm} \calS_\pm$.  Therefore, in order to show that
$\calS_\pm$ is nontorsion, it suffices to produce a generator of
$\chr_\La \calS$ whose image in some $\La_\pm$ is zero.

As in Proposition \ref{ppn-bound}, $C$ is representable by a complex
of the form $[\Phi \stackrel{u}{\hookrightarrow} \Psi]$.  Applying
Proposition \ref{mr-rep}, let us represent $C$ once and for all by a
basic skew-Hermitian complex $C(M,h)^\bullet$.  ($M$ is an
``organizing module'' for the arithmetic of $T$; cf.\ \cite{mr:org}.)
Recall that $\corank_R S_A(K)$ is assumed to be odd, and that
\[
r = \rank_\La M \equiv \rank_R \calS_K = \corank_R S_A(K) \pmod{2}.
\]
As in \S\ref{sub-mr}, take $\Xi$ to be generated by $\iota$ and
$\sigma$, and obtain from \ref{mr-gen} a generator $\calL$ of
$\chr_\La \calS$, together with a homomorphism $\ep \cn \Xi \to
\{\pm1\}$ describing the action of $\Xi$ on $\calL$.  If $\ep(\sigma)
= -1$, then since $\sigma$ acts trivially on $\La_+$ we must have
$\calL \mapsto 0 \in \La_+$, so we are done.  (In this case, we did
not need to assume that $r$ is odd.)  In the case that $\ep(\sigma) =
+1$, we see that $\ep(\sigma\iota) = \ep(\sigma)\,\ep(\iota) =
1\,(-1)^r = -1$, which forces $\calL \mapsto 0 \in \La_-$, since
$\sigma\iota$ acts trivially on $\La_-$.
\end{proof}

\section{Applications}

Here we apply the cohomological theorem to abelian varieties, modular
forms, and Hida families.  In each case, we give hypotheses on the
object in question that guarantee that the Hypotheses (Symp), (Ord),
(Tam), (Irr), and (Zero) of Section \ref{sub-cohom} hold.

\subsection{Abelian varieties}\label{sub-abs}

The following theorem is due to Mazur--Rubin \cite[Theorem
3.1]{mr:growth}.  Although their statement only includes elliptic
curves, their proof applies {\it verbatim} to any abelian variety with
a prime-to-$p$-degree polarization.

\begin{thm}\label{thm-abs}
Let $B/K_0$ be an abelian variety.  Assume the following hypotheses
hold:
\begin{itemize}
\item[(Symp)] $B$ admits a polarization of prime-to-$p$ degree.
\item[(Ord)] At every place $v$ above $p$, $B$ has good, ordinary
  reduction.
\item[(Tam)] At every place $v$ not above $p$ where $B$ has bad
  reduction, $p$ does not divide the order of the component group of
  the special fiber of the N\'eron model of $B/\calO_{K,v}$ (i.e., the
  Tamagawa number of $B$ at $v$).
\item[(Irr)] No $p$-torsion element in $B(K)$ lies in every formal
  group of $B/\calO_{K,v}$, as $v$ ranges over the places of $K$ lying
  above $p$.
\item[(Zero)] Denoting by $\bbF_v$ the residue field of $K$ at $v$,
  one has $B(\bbF_v)[p] = 0$ for every place $v$ of $K$ above $p$.
\end{itemize}
If the classical $p$-power Selmer group $\mathrm{Sel}_{p^\infty}(B/K)$
of $B/K$ has odd rank, then for some choice of sign $\ep = \pm$, one
has $\corank_{\bbZ_p} \mathrm{Sel}_{p^\infty}(B/L) \geq [L:K]$ for all
finite subextensions $L/K$ of $K_\ep$.
\end{thm}

\begin{proof}
We take $R = \bbZ_p$, and $T = T_pB$, the $p$-adic Tate module of $B$.
The set $S$ consists of the places dividing $p\infty$, and where $B$
has bad reduction.

Let $\wh{B}/K$ be the dual abelian variety.  Fix a polarization $\la
\cn B \to \wh{B}$ such that $p \nmid \deg(\la)$.  The composition $T
\stackrel{\la}{\to} T_p\wh{B} \cong T\sto$ (where the second map is
the Weil pairing) is injective, with cokernel of order
$|\deg(\la)|_p\inv = 1$, and hence is an isomorphism.  This
self-duality is symplectic because the Weil pairing is.  Hence $T$
satisfies (Symp).

At a place $v \in \Sigma$, we let $A_v^+$ consist of the image of the
$p$-power torsion of the formal group of $B$ over $\calO_{K,v}$.  The
compatibility of the Weil pairing with Cartier duality of finite flat
group schemes ensures that these local conditions are Lagrangian
(since $B$ is good ordinary at $v$).  Thus $T$ satisfies (Ord).

For (Tam), we point out that the criterion $V^{I_v}/T^{I_v}
\twoheadrightarrow A^{I_v}$ appearing in Remark \ref{rem-DVR-simplify}
is equivalent to $p$ not dividing the order of the component group
appearing in the statement of the theorem, by \cite[Expos\'e IX,
Proposition 11.2]{sga7}.

Because $A^{G_K} = B(K)[p^\infty]$, the hypothesis that $A^{G_K} \to
\bigoplus_{v \in \Sigma} A_v^-$ be injective means precisely that no
global $p$-torsion point comes comes from every $v$'s formal group.
Thus we have (Irr).

The hypothesis (Zero) holds because $(T_v^-/\fkm)^{G_v} =
B(\bbF_v)[p]$ for $v \in \Sigma$.

As we have just verified the hypotheses of Theorem \ref{thm-cohom}, we
deduce the desired growth of $S_A(L)$ in one of the $K_\pm$.  On the
other hand, by \cite[9.6.7.3]{nekovar}, there is a natural injection
$\text{Sel}_{p^\infty}(B/K) \hookrightarrow S_A(L)$ with finite
cokernel, for $L$ any finite extension of $K$.  The desired growth of
$\text{Sel}_{p^\infty}(B/K)$ thus follows from the growth of $S_A(L)$.
\end{proof}

\begin{rem}
The hypothesis (Irr) of the above theorem holds when, in particular,
$B(K)$ has no $p$-torsion.
\end{rem}

\begin{rem}
The $p$-primary part of the Shafarevich--Tate conjecture would imply
that $\rank_\bbZ B(L) = \corank_{\bbZ_p} \text{Sel}_{p^\infty}(B/L)$.
Thus, conjecturally, the Selmer growth guaranteed above actually means
growth of the of Mordell--Weil rank.
\end{rem}

\subsection{Modular forms}\label{sub-mfs}

Let $\calO$ be the ring of integers in a finite extension $F$ of
$\bbQ_p$, and $N$ a positive integer not divisible by $p$.  Suppose we
are given a normalized eigenform $f \in S_k(\Ga_0(N) \cap
\Ga_1(p),\tau_\cycl^i;\calO)$ of {\it even} weight $k = 2k_0 \geq 2$
that is ordinary at $p$, in the sense that there is a $p$-adic unit
root $\alpha_p$ of its $p$th Hecke polynomial.  Choose a lift of $i$
from $\bbZ/(p-1)$ to $\bbZ/(2(p-1))$, noting in passing that $i$ is
even.

Let $T_f$ be the $p$-adic Galois representation associated to $f$ by
Deligne (with the {\it homological} normalization); it is a free
$\calO$-module of rank $2$ equipped with a continuous, linear
$G_{\bbQ,S}$-action, with $S = \{v \text{ dividing } Np\infty\}$.  It
is well-known that because $f$ is ordinary, $T_f|_{G_p}$ is reducible,
admitting a unique decomposition of the form
\begin{eqnarray}\label{eqn-ord-first}
0 \to \calO(\upsilon\inv\chi_\cycl^{k-1}\tau_\cycl^i) \to T_f|_{G_p}
\to \calO(\upsilon) \to 0,
\end{eqnarray}
with $\upsilon$ the unramified character whose value on $\Frob_p$ is
equal to $\alpha_p$.  (For proofs of these facts, see \cite[Theorems
1(1) and 2]{wiles} and a combination of the items \cite[1.5.5 and
1.6.10]{np}.)  Philosophically, our method applies to the Selmer group
related to the central value of an $L$-function, because it uses the
odd sign of a functional equation to produce a trivial zero there.
Under a general recipe, the Selmer group of $T_f$ is related to the
$L$-value $L(f,1)$, which is not the central value when
$k>2$---rather, $L(f,k/2)$ is.  Accordingly, we must twist $T_f$ by a
power of the Tate motive to make it self-dual.

In fact, one has $j_f\cn T_f(\chi_\cycl^{2-k}\tau_\cycl^{-i})
\stackrel{\sim}{\to} T_f\sto$.  Setting
\[
T := T_f(\chi_\cycl^{1-k_0}\tau_\cycl^{-i/2})
\]
(using our chosen lift of $i$ mod $2(p-1)$), we have $j =
j_f(\chi_\cycl^{k_0-1}\tau_\cycl^{i/2}) \cn T \stackrel{\sim}{\to}
T\sto$.  This duality is symplectic because $k$ is even, by
\cite[1.5.5 and 1.6.10]{np}.  Twisting the exact sequence
(\ref{eqn-ord-first}) by $\chi_\cycl^{1-k_0}\tau_\cycl^{-i/2}$, we
obtain the exact sequence
\begin{equation}\label{eqn-ord}
0 \to \calO(\upsilon\inv\chi_\cycl^{k_0}\tau_\cycl^{i/2}) \to T|_{G_p}
\to \calO(\upsilon\chi_\cycl^{1-k_0}\tau_\cycl^{-i/2}) \to 0,
\end{equation}
producing a unique $G_p$-stable $\calO$-direct summand $T_p^+ \subset
T|_{G_p}$ with Galois action through the character
$\upsilon\inv\chi_\cycl^{k_0}\tau_\cycl^{i/2}$.  Moreover, $T_p^-
\cong \upsilon\chi_\cycl^{1-k_0}\tau_\cycl^{-i/2}$.  The uniqueness
implies that this local condition is Lagrangian.

Consider the following hypotheses on $f$:
\begin{description}
\item[(Irred)] $T/\fkm$ has no $G_K$-fixed vectors.
\item[(Min)] The Serre conductor of the residual representation
$T/\fkm$ is $N$.
\item[(Tame)] $\ell^3 \nmid N$ for all primes $\ell$.
\end{description}
For later use, we briefly collect some implications of these hypotheses.

\begin{ppn}\label{ppn-mf-justify}
Let $f$ be a modular form as described above.
\begin{enumerate}
\item If (Irred) and (Min) hold for $f$, then the conductor of $f$ is
  $N$ or $Np$; in particular, $f$ is new at $N$.
\item If (Tame) holds for $f$, then for $\ell \mid N$ the wild inertia
  group $I_\ell^\text{wild} \subset I_\ell$ acts trivially on $T_f$
  (and hence on $T$).
\end{enumerate}
\end{ppn}

\begin{proof}
First, notice that $\chi_\cycl$, and hence $\tau_\cycl$, is unramified
at $\ell$, and therefore $T_f$ and $T$ are isomorphic as
$I_\ell$-modules.

For (1), the condition (Irred) implies that the representation
associated to $f$ is residually irreducible as a $G_\bbQ$-module, and
hence the Serre conductor is computed using $T_f/\fkm$.  The following
inequalities with $\ell \mid N$ are then easy:
\[
\ord_\ell \text{cond}(T_f/\fkm)
\leq \ord_\ell \text{cond}(T_f \otimes_\calO F)
\leq \ord_\ell N.
\]
The hypothesis (Min) requires that the outer terms be equal;
therefore, the middle term equals the final term.

For (2), we must show that $r := \rank_\calO T^{I_\ell^\text{wild}}$
is equal to $2$.  Certainly $0 \leq r \leq 2$, so we must rule out
$r=0$ and $r=1$.  First, notice that if $r=0$, i.e.\
$T^{I_\ell^\text{wild}} = 0$, then $T^{I_\ell} = 0$ as well, and hence
\[
\ord_\ell N = \frac{1}{1}(2-\rank T^{I_\ell}) + \frac{1}{?}(2-\rank
T^{I_\ell^\text{wild}}) + \cdots > 2,
\]
which is a contradiction.  Also, if $r \geq 1$, then
$I_\ell^\text{wild}$ fixes some vector $v \in T$.  Since $\det T =
\chi_\cycl$ is a trivial $I_\ell$-module, the matrix of the
$I_\ell^\text{wild}$-action with respect to a basis extending $v$ has
the form \parbox[b]{0.541in}{\begin{singlespace} $\begin{pmatrix}1 & *
\\ 0 & 1\end{pmatrix}$ \end{singlespace}}.  But $*$ amounts to a
homomorphism from $I_\ell^\text{wild}$ to a pro-$p$ group, and
$I_\ell^\text{wild}$ is a pro-$\ell$-group, so $* = 0$.  Hence, we
must have $r = 2$.
\end{proof}

The following is the main result of this section.

\begin{thm}\label{thm-mfs}
Let $f \in S_{2k_0}(\Ga_0(N) \cap \Ga_1(p),\tau_\cycl^i;\calO)$ be an
ordinary eigenform as above, with $p \nmid N$.  Assume that the
hypotheses (Irred), (Min), and (Tame) hold for $f$.  Let $K$ be a
number field for which, additionally:
\begin{itemize}
\item $K/\bbQ$ is unramified above $Np$.
\item If $i/2+k_0 \equiv 1 \pmod{p-1}$, then
  $\alpha_p^{[\bbF_v:\bbF_p]} \not\equiv 1 \pmod{\fkm}$ for all
  places $v$ above $p$.
\end{itemize}
Suppose $K$ contains a degree $2$ subfield $K_0$, and assume that
$\corank_\calO S_A(K)$ is odd.  Then for some choice of sign $\ep =
\pm$, one has $\corank_\calO S_A(L) \geq [L:K]$ for all finite
subextensions $L/K$ of $K_\ep$.
\end{thm}

\begin{rem}
If $i/2+k_0 \equiv 1 \pmod{p-1}$ then the other lift of $i$ to
$\bbZ/2(p-1)$ has $i/2+k_0 \equiv (p+1)/2 \pmod{p-1}$, so the second
condition on $K$ never rules out both of the two self-dual twists.
\end{rem}

\begin{proof}
The theorem may be deduced from Theorem \ref{thm-cohom} as soon as the
latter's hypotheses are satisfied.  As we have described above, $T$
satisfies (Symp) and (Ord) over $\bbQ$, and hence it also does so over
$K_0$.

The freeness in (Tam) is automatic, because $\calO$ is a DVR.  For the
rest of (Tam), by Remark \ref{rem-DVR-simplify} we must show that
$H^1(I_v,T)_\text{tors}^{\Frob_v=1}$ is trivial.  In fact, since
$K/\bbQ$ is unramified above $N$, for every place $v$ of $K$ lying
over $\ell \in \Sigma'$ we have $I_v = I_\ell$.  Thus, it suffices to
show that $H^1(I_\ell,T)_\text{tors} = 0$.

As in Proposition \ref{ppn-mf-justify}(2), the condition (Tame)
implies that for every prime number $\ell \mid N$, the inertia group
$I_\ell \subset G_\ell$ acts through its tame quotient
$I_\ell/I_\ell^\text{wild}$.  Fix a generator $t \in
I_\ell/I_\ell^\text{wild}$.  The $I_\ell$-module $T$ is
pro-$p$-finite, hence $H^1(I_\ell^\text{wild},M) = 0$, so that
inflation induces the first of the isomorphisms:
\[
H^1(I_\ell,T) \stackrel{\sim}{\leftarrow}
H^1(I_\ell/I_\ell^\text{wild},T) \cong T_{I_\ell}.
\]
Reducing mod $\fkm$, we observe that
\begin{equation}\label{eqn-inertia}
H^1(I_\ell,T)/\fkm \cong T_{I_\ell}/\fkm = T/(t-1,\fkm) =
(T/\fkm)_{I_\ell}.
\end{equation}
Following the conductor identity in Proposition
\ref{ppn-mf-justify}(1),
\begin{equation}\label{eqn-minimal}
\dim_F (T \otimes_\calO F)_{I_\ell} = \dim_k (T/\fkm)_{I_\ell}.
\end{equation}
Moreover, since $f$ is new at $N$, $T$ is ramified at $\ell$, and so
the left hand side of the above equation is $0$ or $1$.

In the case of a $0$, Nakayama's lemma and Equation \ref{eqn-inertia}
show that $H^1(I_\ell,T) = 0$.

Otherwise, in the case of a $1$, we argue as follows.  For any
finitely generated module $M$ over the Noetherian local ring $R$, the
quantity $d = \dim_k M/\fkm$ is the minimal integer for which there
exists a surjection $R^{\oplus d} \twoheadrightarrow M$.  Thus
$T_{I_\ell}$ must be cyclic.  A cyclic module is either free or it is
torsion over some generic component, and the former must hold for
$T_{I_\ell}$.  Otherwise we would have $T_{I_\ell} \otimes_\calO F =
0$, contradicting the $1$ in Equation \ref{eqn-minimal} because $(T
\otimes_\calO F)_{I_\ell} = T_{I_\ell} \otimes_\calO F$.  Thus
$H^1(I_\ell,T)_\text{tors} = 0$.

Putting together the two cases, we see that (Tam) holds in any case.

Now we consider the condition (Irr).  The map $A[\fkm] \to \bigoplus_v
A_v^-$ is nonzero because $A[\fkm] \twoheadrightarrow A_v^-[\fkm]$ for
every $v \in \Sigma$.  Therefore, if the kernel of $A[\fkm] \to
\bigoplus_{v \in \Sigma} A_v^-$ is nonzero, then it must be
$1$-dimensional over $\calO/\fkm$.  But such a subspace would be
$G_{\bbQ,S_0}$-stable, contradicting (Irred).  Since $A[\fkm] \to
\bigoplus_{v \in \Sigma} A_v^-$ is injective, so is its restriction to
$A[\fkm]^{G_{K,S}}$, which gives (Irr).

Finally, we show that (Zero) holds.  Consider the restriction of
$T/\fkm$ to $G_p$.  Since the exact sequence (\ref{eqn-ord}) defines
$T_p^-$, we must show the nontriviality of the character
$\upsilon\chi_\cycl^{1-k_0}\tau_\cycl^{-i/2} \bmod \fkm$ (when further
restricted to $G_v \subseteq G_p$ for $v$ lying over $p$).

The character $\upsilon$ is unramified.  The residual character
$\chi_\cycl^j \bmod \fkm$ is identically equal to $\tau_\cycl^j$, and
on $G_{\bbQ_p}$ it is unramified precisely for $j \equiv 0
\pmod{p-1}$, in which case it is trivial.  So, if $i/2+k_0 \not\equiv
1 \pmod{p-1}$ then our character is always ramified, and hence its
restrictions to the $G_v$ (for $v \in \Sigma$) are nontrivial since
$K/\bbQ$ is unramified at $p$.  Hence, in this case, (Zero) holds.

In the case when $i/2+k_0 \equiv 1 \pmod{p-1}$, our character equals
$\upsilon \bmod \fkm$.  Since $\Frob_p \in G_p/I_p$ acts through
$\upsilon$ via $\alpha_p$, the element $\Frob_v \in G_v/I_v$ acts via
$\alpha_p^{[\bbF_v:\bbF_p]}$.  Thus, by our additional hypothesis,
$\upsilon$ is nontrivial, and (Zero) holds in this case as well.  The
theorem follows.
\end{proof}

\begin{rem}
Let $B = B_f$ be a modular elliptic curve with good ordinary reduction
at $p$.  Then the hypotheses of Theorem \ref{thm-abs} on $B$ are
slightly weaker than the hypotheses of Theorem \ref{thm-mfs} on $f$.
Namely, (Irred) implies (Irr), and if (Irred) holds then (Min) and
(Tame) imply (Tam).  Both reverse implications do not necessarily
hold.  The asymmetry is due to the more complicated nature of
guaranteeing the vanishing of the Tamagawa numbers of a modular form,
relative to the corresponding claim for an abelian variety.
\end{rem}

\subsection{Hida families}\label{sub-hidas}

Fix a positive integer $N$ prime to $p$, an even integer $i$ modulo
$2(p-1)$, and a finite extension of $\bbQ_p$ with ring of integers
$\calO$.  Write $\La$ for the Iwasawa algebra
$\calO[\![\bbZ_p^\times]\!]$, and use brackets $\bra{\ }$ to denote
grouplike elements in it.  (We will not use the $\La$ of
\S\ref{sub-notation} here.)

We take for $R$ a completion of Hida's $\La$-adic Hecke algebra
$\fkh^{\ord}_\infty(\Ga_0(N),\tau_\cycl^i;\calO)$ at a maximal ideal.
It is finite free as a $\La$-module.  Thus, the classical points $\fkp
\in \Spec(R)^\alg$ correspond to $p$-stabilized normalized ordinary
eigenforms
\[
f_\fkp \in S_{k_\fkp}(\Ga_0(N) \cap \Ga_1(p^{c(\fkp)}),
\tau_\cycl^{i-k_\fkp} \psi_\fkp; \calO_\fkp)
\]
lying in a fixed congruence class determined by the maximal ideal.  We
assume that $\Frac \calO$ is algebraically closed in each factor of $R
\otimes_\La \Frac \La$.  (Otherwise, enlarge $\calO$ to ensure this
condition.)

We define $\Theta = \bra{\ga_\cycl^{1/2}}\,\tau_\cycl^{i/2-1}$ (cf.\
\S\ref{sub-notation}).  By knowledge of the diamond operators, the
composite $\mu_{p-1} \subset \bbZ_p^\times \stackrel{\bra{\
}}{\hookrightarrow} \La \to R$ is given by $\mu_{p-1} \ni a \mapsto
a^{i-2} \in \bbZ_p \subseteq R$.  Therefore, $\Theta^2 =
\bra{\ga_\cycl}\,\tau_\cycl^{i-2}$ is the composite of
$\bra{\chi_\cycl}$ with the structure map $\La \to R$; in other words,
$\Theta$ is a square root of $\bra{\chi_\cycl}$ with values in $R$.

Let $T_R$ denote the Galois representation associated to $R$ by Hida
(with the {\it homological} normalization).  It is a free $R$-module
of rank $2$, equipped with a continuous, linear $G_{\bbQ,S}$-action,
where $S$ consists of the places of $\bbQ$ dividing $Np\infty$.  It
has the property that $T_R \otimes_R \calO_\fkp \cong T_{f_\fkp}$ for
all $\fkp \in \Spec(R)^\alg$.  It admits a perfect,
Galois-equivariant, skew-symmetric pairing $T_R \otimes_R T_R \to
R(\bra{\chi_\cycl}\chi_\cycl)$.  Its restriction to $G_p$ is
reducible, sitting in an exact sequence
\begin{equation}\label{eqn-hida-seq-old}
0 \to R(\Upsilon\inv\bra{\chi_\cycl}\chi_\cycl) \to T_R|_{G_p} \to
R(\Upsilon) \to 0,
\end{equation}
where $\Upsilon$ is the unramified $R$-valued character of $G_p$ whose
value on $\Frob_p$ is the Hecke operator $U_p \in R$.  (For details,
see \cite[Chapter 1]{np}.  Strictly speaking, the above is only known
to hold under the assumption of (Irred) below, which will be in force
for our theorem.)

As in the preceding section, the representation $T_R$ corresponds to
the $L$-values $L(f_\fkp,1)$, so we twist it to correspond to the
$L(f_\fkp,k_\fkp/2)$.  We define
\[
T := T_R \otimes_R \Theta\inv.
\]
Twisting the duality pairing $T_R \otimes_R T_R \to
R(\bra{\chi_\cycl}\chi_\cycl)$ by $\Theta^{-2} =
\bra{\chi_\cycl}\inv$, one obtains a duality pairing on $T$ with
values in $R(1)$.  By \cite[1.6.10]{np}, this pairing is symplectic.
Moreover, twisting Equation \ref{eqn-hida-seq-old} by $\Theta\inv$, we
see that the restriction of $T$ to any decomposition group
$G_{\bbQ_p}$ at $p$ sits in an exact sequence
\[
0 \to R(\Upsilon\inv\,\Theta\,\chi_\cycl) \to T|_{G_{\bbQ_p}} \to
R(\Upsilon\,\Theta\inv) \to 0.
\]
It is easy to see that this exact sequence is self-dual under the
pairing just mentioned.  We define $T_p^+$ to be the image of
$R(\Upsilon\inv\,\Theta\,\chi_\cycl)$ in $T$.

Here are our hypotheses:
\begin{description}
\item[(Irred)] $T/\fkm$ has no $G_K$-fixed vectors.
\item[(Min)] The Serre conductor of the residual representation
$T/\fkm$ is $N$.
\item[(Tame)] $\ell^3 \nmid N$ for all primes $\ell$.
\end{description}

Our theorem closely resembles the corresponding result, Theorem
\ref{thm-mfs}, for individual modular forms.

\begin{thm}\label{thm-hidas}
Let $R$ be a completion of the Hida--Hecke algebra of tame level
$\Ga_0(N)$ and character $\tau_\cycl^i$, and let $T$ be the associated
twisted representation, as above.  Assume that (Irred), (Min), and
(Tame) hold.  Let $K$ be a number field for which, additionally:
\begin{itemize}
\item $K/\bbQ$ is unramified above $Np$, and for each place $v$ of $K$
  lying over $N$, we have $U_\ell^{[\bbF_v:\bbF_\ell]} \notin 1 +
  \fkm$.
\item If $i/2 \equiv 1 \pmod{p-1}$, then $U_p^{[\bbF_v:\bbF_p]} \notin
  1 + \fkm$ for all places $v$ lying over $p$.
\end{itemize}
Suppose $K$ contains a degree $2$ subfield $K_0$, and assume that
$\corank_R S_A(K)$ is odd.  Then for some choice of sign $\ep = \pm$,
one has $\corank_R S_A(L) \geq [L:K]$ for all finite subextensions
$L/K$ of $K_\ep$.
\end{thm}

\begin{rem}
If $i/2 \equiv 1 \pmod{p-1}$ then the other lift of $i$ to
$\bbZ/2(p-1)$ has $i/2 \equiv (p+1)/2 \pmod{p-1}$, so the second
condition on $K$ never rules out both of the two self-dual twists.
\end{rem}

\begin{proof}
We will deduce our desired claim from Theorem \ref{thm-cohom}.  To
begin with, the ring $R$ is a complete Noetherian local ring; by
(Irred) and \cite[1.5.2--1.5.4]{np}, $R$ is Gorenstein, and the
``residual representation'' of the Hida family is realized by
$T/\fkm$.

As described above, $T$ satisfies (Symp) and (Ord) over $\bbQ$, and,
consequently, also over any finite extension of $\bbQ$.

We now consider (Tam).  As in the case of modular forms, (Tame)
implies that for $v$ lying over $\ell \in \Sigma'$, the inertia group
$I_\ell$ acts on $T$ through its tame quotient, and
\[
H^1(I_\ell,T) \cong T/(t-1) \qquad \text{and} \qquad
H^1(I_\ell,T/\fkm) \cong T/(\fkm,t-1),
\]
 for $t \in I_\ell/I_\ell^\text{wild}$ a generator.  Then the same
reasoning as in the proof of Theorem \ref{thm-mfs} using (Irr) and
(Min) shows that $H^1(I_\ell,T)$ is either zero or free of rank one.

Suppose that $H^1(I_\ell,T)$ is nonzero.  Note that under our
hypothesis (Min), in the notation of \cite{epw}, one has
$\bfT_N^\text{new} = \bfT_N$.  Therefore, \cite[Lemma 2.6.2]{epw}
implies that $\text{Frob}_\ell$ acts on the line $(T \otimes_R \Frac
\La)_{I_\ell}$ through the Hecke operator $U_\ell$; the same is
therefore true on the lattice $T_{I_\ell}$.  Since $v$ is unramified
over $\ell$ one has $I_v = I_\ell$, and the preceding calculation
shows that $\Frob_v-1$ acts bijectively on $H^1(I_v,T)$ precisely when
$U_\ell^{[\bbF_v:\bbF_\ell]} \notin 1 + \fkm$, which was our
additional hypothesis.

For the freeness part of (Tam), recall that $T$ is self-dual, so that
$T \cong T\sto$.  We compute from this (noting that $I_\ell$ acts
trivially on $R(1)$) that
\begin{eqnarray}\label{eqn-hida-invariants}
T^{I_\ell}
 & \cong & \Hom(T,R(1))^{I_\ell}
 = \ker\left(t-1 \mid \Hom(T,R(1))\right) \nonumber \\
 & = & \displaystyle
       \Big\{f \cn T \to R(1) \ \Big| \ 0 = [(t-1)f](T) =
        f((t\inv-1)T) \Big\} \\
 & = & \Hom(T/(t\inv-1),R(1)) = \Hom(T_{I_\ell},R(1)). \nonumber
\end{eqnarray}
Since $T_{I_\ell}$ is free in all cases, we deduce from the above
identity that $T^{I_\ell}$ is free too.  Hence (Tam) holds.

The verification of the hypotheses (Irr) and (Zero) is accomplished in
the exact same manner as in the proof of Theorem \ref{thm-mfs}, so we
omit it.  The only (cosmetic) difference is that the $G_p$-action on
$T_v^-/\fkm$ here is through $\Upsilon\tau_\cycl^{1-i/2} \bmod \fkm$,
which plays the same role as $\upsilon\tau_\cycl^{1-i/2-k_0} \bmod
\fkm$ in the proof of Theorem \ref{thm-mfs}.
\end{proof}

Finally, we mention the following result that explicates the
relationship between the theorems for modular forms and their Hida
families.

\begin{ppn}
Assume $R$ is a UFD (e.g, it is regular).  Assume that, for all places
$v$ of $K$ above $N$, we have $U_\ell^{[\bbF_v:\bbF_\ell]} \notin 1 +
\fkm$.  Then the remaining hypotheses of Theorem \ref{thm-hidas} are
satisfied by $R$ if and only if the hypotheses of Theorem
\ref{thm-mfs} are satisfied by any classical, even weight modular form
$f_\fkp$ in the family with $c(\fkp) \leq 1$.  In this case, the one
has $\calS(f_\fkp) = (\calS(R) \otimes_R R/\fkp) \otimes_{R/\fkp}
\calO_\fkp$, with the obvious notation.
\end{ppn}

\begin{proof}
Using the description of $\calS$ provided by Proposition
\ref{ppn-bound}, this is an easy consequence of Nekov\'a\v{r}'s
control theorem \cite[8.10.10]{nekovar} and Nakayama's lemma, since
every height-one prime $\fkp \in \Spec R$ is principal.  (Cf.\ the
comments at the end of \S\ref{sub-global}, where Equation
\ref{eqn-hyp-control} holds because of Equation
\ref{eqn-hida-invariants} above, using the fact that $T_{I_\ell}$, if
nonzero, is free of rank one.)  As in the verification of (Zero) in
the above proof, note the slight renormalization of the relevant power
of $\tau_\cycl$.
\end{proof}

It follows from the above proposition that, when its conditions are
satisfied, for each sign $\ep = \pm$ one has
\[
\rank_{\La(R)_\ep} \calS(R)_\ep = \rank_{\La(f_\fkp)_\ep}
\calS(f_\fkp)_\ep
\]
for all but finitely many classical, even-weight $\fkp$, and this rank
is positive for at least one sign.  Thus, the Selmer growth of modular
forms, as produced in this paper, occurs ``uniformly in the Hida
family'', with finitely many exceptions.

\section{Examples}

Here we show concrete instances of the applications of the preceding
chapter.

\subsection{Applications of the Cassels--Tate pairing}\label{sub-CT}

We remind the reader that we had no need for \Nekovar's higher
Cassels--Tate pairing apparatus in the preceding chapters; this is
still the case for the examples of modular abelian varieties below.
However, we admit applications of the Cassels--Tate pairing below to
produce specific numerical examples of modular forms of higher weight,
and of Hida families.  Precisely, we employ the parity conjecture as
proved in \cite[12.2.6]{nekovar} by means of Cassels--Tate in order to
detect odd Selmer parity; and, in order to extend the computations to
Hida families, we refer to a direct application of Cassels--Tate
\cite[10.7.5.ii]{nekovar}.  (The latter result, on its own, is enough
to prove a fine Selmer growth theorem.)

Assume we are given an even-weight eigenform $f \in S_k(\Ga_0(N) \cap
\Ga_1(p), \tau_\cycl^i; \calO)$ that is ordinary at $p$ with $p \nmid
N$ as in \S\ref{sub-mfs}, and an imaginary quadratic field $K/\bbQ$
that is unramified at $Np$.  Assume all the hypotheses of Theorem
\ref{thm-mfs} hold except possibly the parity condition: namely,
(Irred), (Min), and (Tame) hold, and that if $i/2+k/2 \equiv 1
\pmod{p-1}$ then $\alpha_p^{[\bbF_\fkp:\bbF_p]} \not\equiv 1
\pmod{\fkm}$ (a choice of $i/2 \in \bbZ/(p-1)$ satisfying this can
always be made).  The question is for which choices of $K$ one has
$\corank_\bbZp S_{A_f}(K)$ odd.

To answer this, we first invoke the parity conjecture
\cite[12.2.6]{nekovar}, whose proof makes use of Cassels--Tate.  This
theorem implies, in particular, that
\[
\corank_\calO S_{A_f}(L) \equiv \ord_{s=k/2} L(f/L,s) \pmod{2}
\]
for $L = \bbQ \text{ or } K$.  In the case $L=K$, denoting by $\chi_K$
the quadratic Dirichlet character corresponding to $K$, one has
\[
L(f/K,s) = L(f/\bbQ,s)L(f_{\chi_K}/\bbQ,s).
\]
Hence, the Selmer parity is the composite of the analytic parities of
$f$ and $f_{\chi_K}$.  Moreover, if $r$ is the analytic rank of $f$
over $L$, it is well-known that $(-1)^r$ is the sign of the functional
equation of $L(f/L,s)$, and when that $L=\bbQ$ this sign is equal to
$(-1)^{k/2}w_{N(f)}(f)$, where $w_{N(f)}(f) = \pm1$ is the eigenvalue
of the Atkin--Lehner $W_{N(f)}$-operator acting on $f$.  (We stress
the dependence of $N$ on $f$ to avoid confusion when applying these
remarks to $f_{\chi_K}$.)  Thus, writing $w(f,K) =
w_{N(f_{\chi_K})}(f_{\chi_K})/w_{N(f)}(f)$ for the {\it change} in
sign caused by twisting $f$ by $\chi_K$, the sign of $L(f/K,s)$ is
computed to be
\[
w_{N(f)}(f)\,w_{N(f_{\chi_K})}(f_{\chi_K}) =
w_{N(f)}(f)\,w(f,K)w_{N(f)}(f) = w(f,K),
\]
as $w_{N(f)}(f) = \pm1$.  But $w(f,K)$ is easily computed by
\cite[3.63(2) and 3.65]{shimura} to be $\chi_K(-N(f))$.  As $K$ is
imaginary quadratic, we find odd Selmer corank precisely when
$\chi_K(N(f)) = 1$.

As for the Hida family passing through $f$, let $R$ be the
corresponding completion of
$\fkh_\infty^{\ord}(\Ga_0(N),\tau_\cycl^i;\calO)$, and assume moreover
that for $U_\ell^{[\bbF_v:\bbF_\ell]} \notin 1+\fkm_R$ for each place
$v$ of $K$ lying over $\ell$ dividing $N$.  We want to find quadratic
fields where odd Selmer corank occurs.  By \Nekovar's Cassels--Tate
apparatus (see \cite[10.7.15.ii]{nekovar}), the $R$-corank of
$S_{A_R}(K)$ is equal to the $\calO$-corank $S_{A_f}(K)$.  In
particular, there is no further restriction on $K$.

\subsection{$\Delta$ at $p=11$}\label{delta11}

Let $\Delta \in S_{12}(\SL_2(\bbZ);\bbZ)$ be the discriminant form,
with
\[\begin{split}
\Delta(q) = \sum_{n \geq 1} \tau(n)q^n
& =  q - 24q^2 + 252q^3 - 1472q^4 + 4830q^5 - 6048q^6 - 16744q^7 \\
& \quad + 84480q^8 - 113643q^9 - 115920q^{10} + 534612q^{11} - \cdots. \\
\end{split}\]
It is the cusp eigenform of lowest weight with level $1$, and $p=11$
is its first ordinary prime.

We take $K_0 = \bbQ$, and ask which imaginary quadratic fields $K$
satisfies the hypotheses of Theorem \ref{thm-mfs}.  In the notation of
\S\ref{sub-mfs}, one has $(k,N,i) = (12,1,0)$.  Since $N=1$, there is
no need to worry about (Min) and (Tame).  For (Irred), we need to know
for which imaginary quadratic fields $K$ the action of $G_K$ on
$T/\fkm$ leaves no vectors fixed.  For this, we recall that, by the
calculations of Serre--Swinnerton-Dyer \cite{serre,swin}, the action
of $G_\bbQ$ on $T_f/\fkm$ gives a {\it surjection} $G_\bbQ \to
\GL_2(\bbF_{11})$.  Invoking the following elementary result, we
conclude that $\Delta$ satisfies (Irred) for all $K$.

\begin{ppn}\label{ppn-gp-thy}
Let $q$ be a prime power, $G$ a group, and $\rho \cn G \to
\GL_2(\bbF_q)$ a homomorphism with image $G_0$.  Suppose that
$[\GL_2(\bbF_q):G_0] < (q+1)/2$.  Then for no index $2$ subgroup $H
\subset G$ and no character $\chi \cn G \to \bbF_q^\times$ is there a
fixed vector of $\bbF_q^{\oplus 2}$ under $(\rho \otimes \chi)(H)$.
\end{ppn}

\begin{proof}
Such a fixed vector would result in $\rho(H) \subseteq G_1 :=
\parbox[b]{0.542in}{\begin{singlespace} $\begin{pmatrix}* & * \\ 0 &
*\end{pmatrix}$ \end{singlespace}}$, and the latter group has
$q(q-1)^2$ elements.  But
\[
\#G_0 = \frac{\#\GL_2(\bbF_q)}{[\GL_2(\bbF_q):G_0]} > 
\frac{(q^2-1)(q^2-q)}{(q+1)/2} = 2q(q-1)^2 = 2\#G_1,
\]
whence $[G_0:G_1] > 2 = [G:H]$, a contradiction.
\end{proof}

Since $i = 0 \in \bbZ/10$, one has $i/2 = 0 \text{ or }5 \in \bbZ/10$.
Its $11$-Hecke polynomial, modulo $11$, is $X^2-X$, so its unit root
$\alpha_{11} \in \bbZ_{11}$ satisfies $\alpha_{11} \equiv 1
\pmod{11}$.  This forces us to take the self-dual twist of $T_\Delta$
with $i/2 = 0$ (with $i/2+k/2=6$) and not $i/2 = 5$ (with $i/2+k/2 =
1$).  The other restriction listed in Theorem \ref{thm-mfs} is that
$K/\bbQ$ be unramified above $p = 11$.  By the discussion of
\S\ref{sub-CT}, we are now limited to those quadratic fields $K$ with
$\chi_K(N(\Delta)) = 1$.  But $N(\Delta) = 1$, so this is no
restriction.

Finally, since $K/\bbQ$ is abelian, Kato's Euler system (see
\cite{kato}) implies that $\calS_+$ is $\La_+$-torsion.  This
eliminates the possibility of growth with $\ep=+$.  We summarize:

\begin{thm}
Let $K/\bbQ$ be any imaginary quadratic field whose discriminant is
prime to $11$.  Let $A = T_\Delta(\chi_\cycl^{-5}) \otimes_{\bbZ_{11}}
\bbQ_{11}/\bbZ_{11}$.  Then $S_A(L)$ has corank bounded below by
$[L:K]$ for $L$ ranging through the finite subextensions of the
anticyclotomic $\bbZ_{11}$-extension of $K$.
\end{thm}

Most of the above work carries over with little change for Hida
theory.  It is well-known that $\La \stackrel{\sim}{\to}
\fkh_\infty^{\ord}(\SL_2(\bbZ);\bbZ_{11})$, i.e.\ the Hida family
passing through $\Delta$ is as small as possible.  (This can be
accomplished, for example, by enumerating cusp eigenforms of
appropriate weight, level and character congruent modulo $11$ to
$\Delta$, and finding only one in each weight, level and character.)

We seek imaginary quadratic fields $K$ to which Theorem
\ref{thm-hidas} applies.  Since $(N,i) = (1,12)$, we still need not
worry about (Tam) and (Min).  Since $T/\fkm$ for the Hida family is
identical to that for $f$, we still see that every quadratic field
satisfies (Irred).  The theorem requires that $K$ be unramified at
$11$.  We easily treat the numerology: $i/2 = 1 \text{ or } 6 \in
\bbZ/10$, and we must take $i/2 = 6$, because $U_p \equiv \alpha_p
\equiv 1 \pmod{\fkm}$.  By \S\ref{sub-CT}, the same quadratic fields
$K$ are valid as those for $\Delta$ itself.  Finally, we point out
that Kato's result implies that $\calS_+$ is
$\Lambda(\text{Hida})_+$-torsion, because $\calS(\text{Hida})_+
\twoheadrightarrow \calS(\Delta)_+$, and the latter is
$\Lambda(\Delta)_+$-torsion, in the obvious notation.  We summarize:

\begin{thm}
Let $K/\bbQ$ be any imaginary quadratic field whose discriminant is
not divisible by $11$.  Let $A =
T_\text{Hida}(\bra{\ga_\cycl^{-1/2}}\tau_\cycl^{-5}) \otimes_{\bbZ_p}
\bbQ_p/\bbZ_p$, where $T_\text{Hida}$ is the Galois representation
constructed by Hida with values in his $\La$-adic Hecke algebra $R$
with $(p,N,i) = (11,1,12)$.  Then $S_A(L)$ has corank bounded below by
$[L:K]$ for $L$ ranging through the finite subextensions of the
anticyclotomic $\bbZ_p$-extension of $K$.
\end{thm}

\subsection{$y^2+y=x^3-x$ over $\bbQ(\sqrt{-3})$}

The following example is taken from \cite[Example 10.10]{mr:org}, and
hence is well-known.  Let $E$ be the elliptic curve over $\bbQ$
defined by
\[
y^2+y=x^3-x.
\]
Instead of fixing $p$ and varying $K$, we fix $K = \bbQ(\sqrt{-3})$
and ask when Theorem \ref{thm-abs} applies.  It is claimed in
\cite{mr:org} that $E(K) = E(\bbQ) = \bbZ \cdot (0,0)$ and
$\Sha(E/K)=0$; in particular, $\text{Sel}_p(E/K) \approx
\bbQ_p/\bbZ_p$ has $\bbZ_p$-corank one for any prime $p$.  Let
$\sum_{n \geq 1} a_n q^n$ be the $q$-expansion of the normalized
newform $f$ of level $\Ga_0(37)$ corresponding to $E$.

For (Symp), all elliptic curves are principally polarized.  We have
(Ord) when p>3 and $a_p \neq 0$ (since $p=2,3$ are supersingular for
$E$).  For the Tamagawa obstruction, we note that the residual
representation of $f$ at any $p \neq 37$ has Serre conductor $37$
(because otherwise level-lowering would produce a weight-$2$ cusp form
for $\SL_2(\bbZ)$), which is enough to force
$H^1(I_{37},T)_\text{tors}=0$.  Moreover, $37$ is unramified in $K$,
so this calculation applies to the inertia groups of $K$ lying over
$37$.  Since $E(K) \approx \bbZ$, (Irr) holds for all $p$.  Finally,
(Zero) holds over $\bbQ$ when $a_p \neq 1$, by the Weil bound, and
hence over $K$ when either $a_p$ or $a_p^2 \neq 1$, according to
whether $p$ is split or inert in $K$, respectively.

Once again, by Kato, it is known that $\calS_+$ is $\La_+$-torsion.
We summarize:

\begin{thm}
Let $E$ and $K$ be as given above.  Then, for any prime $p>3$ with $p
\neq 37$ and $a_p \neq 0,\pm 1$, the $\bbZ_p$-coranks of
$\text{Sel}_{p^\infty}(E/L)$ are bounded below by $[L:K]$ for $L$
ranging through the finite subextensions of the $p$-anticyclotomic
extension of $K$.
\end{thm}

\subsection{The unique newform in $S_4(\Ga_0(5))$}

Let $f \in S_4(\Ga_0(5);\bbZ)$ be the unique newform, which has
\[\begin{split}
f(q) = \sum_{n \geq 1} a_nq^n
& = q - 4q^2 + 2q^3 + 8q^4 - 5q^5 - 8q^6 + 6q^7 - 23q^9 + 20q^{10}
 + 32q^{11} \\
& \quad + 16q^{12} - 38q^{13} - 24q^{14} - 10q^{15} - 64q^{16} + 26q^{17}
 + 92q^{18} + \cdots.
\end{split}\]
It is the cusp eigenform of lowest level with weight $k=4$.  It is
clear from the $q$-expansion above that all primes $p \leq 17$ except
$p=2,5$ are good ordinary for $f$; in fact, $p=2,5$ are the only
nonordinary primes for $f$ that are less than $1000$.

Suppose $f$ is ordinary at $p$ (so in particular $p \neq 2,5$).  Then
$f$ satisfies (Min) because there are no level $1$ forms to
level-lower to, and (Tam) because $N = 5$ is prime.  By Proposition
\ref{ppn-gp-thy}, we have (Irred) for all quadratic fields $K$,
granted $p>6$ and $p$ is not exceptional for $f$ (i.e., since the
weight of $f$ is $4$, that the mod $p$ representation associated to
$f$ surjects onto the matrices in $\GL_2(\bbF_p)$ whose determinants
are cubes in $\bbF_p$).  Ribet's generalization \cite{ribet} of
Serre--Swinnerton-Dyer's arguments \cite{serre,swin} proves that there
are only finitely many exceptional $p$; stepping through his paper in
this case, we computed that there are no exceptional $p$ with $p >
19$.  For $p>2$, one certainly has $2,2+(p-1)/2 \not\equiv 1
\pmod{p-1}$, so both self-dual twists work.  Once again, Kato's work
precludes Selmer growth in the cyclotomic extension.  We summarize:

\begin{thm}
let $f$ be as given above.  Let $p$ be an ordinary prime for $f$ with
$p > 19$, and let $T$ be either of the two self-dual twists of the
$p$-adic Galois representation $T_f$ associated to $f$.  Let $K$ be
any imaginary quadratic field in which $5$ and $p$ are unramified,
satisfying $\chi_K(5) = 1$.  Then for all finite subextensions $L/K$
within the anticyclotomic $\bbZ_p$-extension of $K$, one has
$\corank_\bbZp S_A(L) \geq [L:K]$.
\end{thm}

\appendix

\section[Selmer complexes]{Appendix: Selmer complexes}

In this appendix we give a crash course in Selmer complexes, with an
aim towards explaining how the complex $C$ of Proposition
\ref{ppn-bound} is constructed.  This construction is not necessary
for our results, and is provided only for the convenience of the
curious reader; therefore, for the sake of brevity, we sketch the
constructions of our objects, and those facts that will be of use to
us, only in as much generality as is needed.  For much more
information, including proofs of our claims, further properties, and
applications, the reader is referred to \Nekovar's original work
\cite{nekovar}.\footnote{And a hearty referral it is!  The book is
carefully written and includes copious useful details.}

\subsection{Additional notation}

Throughout this appendix, \S\ref{sub-notation} is in force.

There are some full subcategories of $R\ddmod$ that are of interest to
us.  We write $R\ddfl$ (resp.\ $R\ddft$, $R\ddcoft$) for the category
consisting of finite-length (resp.\ Noetherian, Artinian) $R$-modules.
Since $R$ is Noetherian, it is customary to call Noetherian
$R$-modules ``of finite type''; we will adopt the convention of
calling Artinian $R$-modules ``of cofinite type'', thus explaining the
notation.  Objects in the three categories of this paragraph are
naturally endowed with linear topologies.

If $G$ is a group, then $R[G]\ddmod$ is the category of (linearly
topologized) $R$-modules equipped with continuous, linear $G$-actions.
For $\star$ one of the conditions $\ddfl,\ddft,\ddcoft$, we write
$R[G]\ddmod_{R\star}$ for the full subcategory whose underlying
$R$-modules belong to $R\star$.

Throughout this section, we use the letter $M$ (resp.\ $T$, $A$)
exclusively to denote an element of $R\ddfl$ (resp.\ $R\ddft$,
$R\ddcoft$), or perhaps $R[G]\ddmod_{R\ddfl}$ (resp.\
$R[G]\ddmod_{R\ddft}$, $R[G]\ddmod_{R\ddcoft}$).  The letter $X$
stands for any of $M,T,A$, unless otherwise specified.  Written with a
big dot, $X^\bullet$ means either $M^\bullet$, $T^\bullet$, or
$A^\bullet$, with each component satisfying the appropriate finiteness
condition, and is assumed to be a {\it bounded} complex.

\subsection{Continuous cochain complexes}

Fix a profinite group $G$ satisfying the following {\it finiteness
hypothesis}: For every open normal subgroup $H \subseteq G$, all the
$k$-vector spaces $H^*(H,k)$ are finite-dimensional, and moreover $e
:= \text{c.d.}_p(G) < \infty$ (in the sense of cohomology of discrete
modules).  The results that follow will generally only be valid for
$G$ of this type.

We write $C^\bullet(G,M)$ for the {\it
continuous cochain complex} of $G$ with coefficients in $M$, as in
\cite[3.4.1.1]{nekovar}.  Because $M$ is topologically discrete, here
``continuous'' simply means ``locally constant on $G$''.  We put
\[
C^\bullet(G,T) = \llim_n C^\bullet(G,T/\fkm^n) \qquad \text{and}
\qquad C^\bullet(G,A) = \rlim_n C^\bullet(G,A[\fkm^n]).
\]
One extends this notion to complexes by defining the cochain complex
``$C^\bullet(G,X^\bullet)$'' of $X^\bullet$ to be the total complex of
the bicomplex ``$C^\bullet(G,X^\bullet)$''.  (We will never have need
to refer to the latter, so there should be no confusion of notation.)

The association
\[
X^\bullet \mapsto \bfR\Ga(G,X^\bullet)
  := [C^\bullet(G,X^\bullet)] \in \bfD(R)
\]
factors uniquely through an exact functor $\bfD(R[G]\ddmod_{R\ddfl})
\to \bfD(R)$ when $X=M$ (resp.\ $\bfD(R[G]\ddmod_{R\ddft}) \to
\bfD(R)$ when $X=T$, $\bfD(R[G]\ddmod_{R\ddcoft}) \to \bfD(R)$ when
$X=A$).  Passing to cohomology groups, we get the usual ``continuous
hypercohomology'' cohomological $\delta$-functor, which we denote by
$H^*(G,\cdot)$.  (In \cite{nekovar}, the notation $C_\text{cont}$ is
used, but since this is the only type of cohomology we will encounter,
there is no need to continuously remind the reader of any
distinction.)

If $H \subseteq G$ is a closed normal subgroup, then we have an
``inflation'' map
\[
\text{inf}\cn C^\bullet(G/H,X^H) \longrightarrow C^\bullet(G,X),
\]
which is a quasi-isomorphism if $H$ has profinite order not divisible
by $p$.  This leads one to expect the existence of a
``Hochschild--Serre'' morphism of complexes
\[
\text{``}C^\bullet(G/H,C^\bullet(H,X)) \to C^\bullet(G,X)\text{''}.
\]
however, the actions of elements $g,h \in G/H$ on $C^\bullet(H,X)$ are
each only defined {\it up to homotopy}, and the associativity of the
action of $gh$ with those of $g,h$ has {\it not} been verified.  Thus,
the complex $C^\bullet(G/H,C^\bullet(H,X))$ is not yet known to be
well-defined.

The preceding concern has consequences for Iwasawa theory.  We would
like to define, for $H \subseteq G$ a closed normal subgroup with
$G/H$ abelian, complexes
\begin{gather*}
C^\bullet_\text{Iw}(G,H,T)
  := \llim_{\alpha,\text{ cor}} C^\bullet(G_\alpha,T)
\qquad \text{and} \\
C^\bullet(H,A)
  := \rlim_{\alpha,\text{ res}} C^\bullet(G_\alpha,A),
\end{gather*}
the $G_\alpha$ ranging over all open normal subgroups of $G$
containing $H$, representing objects of $\bfD(\La)$ with $\La =
R[\![G/H]\!])$.  However, the action of $G/H$ has not been checked to
make sense, so the right hand sides only {\it a priori} live in
$\bfD(R)$.  To circumvent this, \Nekovar\ has engineered a variant of
Shapiro's lemma that shows that
\begin{gather*}
\llim_{\alpha,\text{ cor}} C^\bullet(G_\alpha,T)
  \cong C^\bullet(G,T \otimes_R \La)
\qquad \text{and} \\
\rlim_{\alpha,\text{ res}} C^\bullet(G_\alpha,A),
  \cong C^\bullet(G,\Hom_R(\La,A)).
\end{gather*}
We thus represent $\bfR\Ga_\text{Iw}(G,H,T)$ and $\bfR\Ga(H,A)$ by the
respective right hand sides of the two above equations, which
ostensibly are equipped with continuous $G/H$-actions, and consider
them as objects of $\bfD(\La)$.  Note that
\begin{gather*}
T \otimes_R \La \in \La[G]\ddmod_{\La\ddft} \qquad \text{and} \\
\Hom_R(\La,A) \in \La[G]\ddmod_{\La\ddcoft},
\end{gather*}
so the above applications of $C^\bullet(G,\cdot)$ make sense.

Under our finiteness hypothesis on $G$, using a Mittag--Leffler
argument, we find that taking $\llim_{\alpha,\text{ cor}}$ is exact,
so no $\bfR\!\llim$ is needed.  Moreover, forming $\llim_\alpha$ and
$\rlim_\alpha$ commute with passing to cohomology groups.  A ``way-out
functors'' argument shows that, for each $i \geq 0$,
\begin{eqnarray*}
H^i(G,M) & \in & R\ddfl, \\
H^i(G,T) & \in & R\ddft, \text{ and} \\
H^i(G,A) & \in & R\ddcoft,
\end{eqnarray*}
and all these groups vanish for $i > e$.  Furthermore, if $T^\bullet
\in \bfD_\text{perf}^{[a,b]}(R)$ (upon forgetting the $G$-action),
then $\bfR\Ga(G,T) \in \bfD_\text{perf}^{[a,b+e]}(R)$.  This applies
particularly to when $T$ is free over $R$, so that $\bfR\Ga(G,T) \in
\bfD_\text{perf}^{[0,e]}(R)$.

\subsection{Intermission: autoduality theorems for $R$}\label{sub-inter}

Since $R$ is a complete Noetherian local ring, Matlis duality theory
takes a simple form.  Let $I_R$ be an injective hull of $k$ in
$R\ddmod$.  Then the rule $X \mapsto D(X) := \Hom_R(X,I_R)$ describes
an exact anti-equivalence of the category $R\ddfl$ with it self, with
$\text{Id}_{R\ddfl} \stackrel{\sim}{\to} D \circ D$.  Hence, taking
limits, $D$ furnishes an exact anti-equivalence $R\ddft
\longleftrightarrow R\ddcoft$.

The finiteness of the residue field $k$ implies that $D$ actually
coincides with the Pontryagin duality functor $X \mapsto
\Hom_{\bbZ_p}(X,\bbQ_p/\bbZ_p)$; this shows that the underlying
$\bbZ_p$-module of $D(X)$ does not depend on $R$.  Moreover, the
objects of $R\ddfl$ are the precisely finite modules, those of
$R\ddft$ are compact, and those of $R\ddcoft$ are discrete.

We required that $R$ be {\it Gorenstein} essentially to guarantee that
Grothendieck duality over $R$ be realized via the complex consisting
of $R$ concentrated in degree $0$.  So, the functor $T^\bullet \mapsto
\scrD(T^\bullet) := \bfR\!\Hom(T,R)$ provides an exact
anti-equivalence of categories $\bfD(R\ddft) \to \bfD(R\ddft)$.  In
the case where $T^\bullet$ is a complex of {\it free} modules, we may
represent $\bfR\!\Hom(T^\bullet,R)$ by the object $\Hom(T^\bullet,R)$
itself.  Thus, $\scrD(T^\bullet)$ exchanges
$\bfD_\text{perf}^{[a,b]}(R)$ with $\bfD_\text{perf}^{[-b,-a]}(R)$.
In this case, including notably when $T^\bullet$ is a single free
module in degree $0$, we will always make this identification,
essentially obliterating the use of derived categories.

Although we will not make use of it, we note that the functors $D$ and
$\scrD$ are related to each other via Grothendieck's local duality
theorem.

Finally, if $\Ga \approx \bbZ_p^n$, then $\La := R\bbGa$ also
satisfies the hypotheses we have made on $R$, and we denote its Matlis
(resp.\ Grothendieck) duality functors by $D_\La$ (resp.\
$\scrD_\La$).

\subsection{Local duality and Selmer complexes}

The arithmetic duality theorems of Poitou--Tate are stated for Galois
modules of finite cardinality.  We present here some results of
\Nekovar's process of ``bootstrapping'' them up to modules over a
larger ring.

Straying somewhat from \S\ref{sub-notation}, we let $K$ be any number
field (forgetting $K_0$), and consider a finite set of places $S$,
partitioned as in \S\ref{sub-notation}, with still $p > 2$.  (The set
of places $S$ of \S\ref{sub-notation} is replaced here by the set of
places of $K$ lying over $S$.)

Let $v$ be a place of $K$.  We may summarize the local duality theorem
as saying that $G_v$ satisfies the finiteness hypothesis with
$\text{c.d.}_p(G_v) = 2$, and that the ``invariant map'' applied to
the cup product has adjoint morphisms
\begin{eqnarray*}
\bfR\Ga(G_v,D(X^\bullet)(1))
  & \to & D\left(\bfR\Ga(G_v,X^\bullet)\right)[-2] \quad \text{and} \\
\bfR\Ga(G_v,\scrD(T^\bullet)(1))
  & \to & \scrD\left(\bfR\Ga(G_v,T^\bullet)\right)[-2]
\end{eqnarray*}
that are isomorphisms in $\bfD(R)$, for $X^\bullet \in
\bfD(R[G_v]\ddmod_{R\star})$ with $\star$ one of $\ddfl$, $\ddft$,
$\ddcoft$, and for $T^\bullet \in \bfD(R[G_v]\ddmod_{R\ddft})$.

The global duality theorem is somewhat more complicated.  The
cohomology over $G_K$ is uncontrollable, so one uses $G_{K,S}$.  And
the cohomology over $G_{K,S}$ is not self-dual, but instead a version
of the cohomology that has been modified by local conditions (i.e.\
the Selmer complex) is self-dual.

A {\it system of local conditions} $\Delta$ for $X^\bullet$ is the
specification of, for each place $v \in S_f$, a complex $U_v^+$ of
$R$-modules and a chain map
\[
i_v^+ \cn U_v^+ \to C^\bullet(G_v,X^\bullet).
\]
The two simplest examples are $\Delta_\emptyset$ (``no conditions'' or
``empty conditions'') with all $i_v^+$ equal to the identity map, and
$\Delta_c$ (``compact support'', or ``full conditions'') with all
$i_v^+$ equal to the map from the zero complex.

For any particular place $v \in S_f$, suppose we are given a complex
of $R[G_v]$-modules $X_v^+ = (X_v^+)^\bullet$ and a $G_v$-morphism
$j_v^+ \cn (X_v^+)^\bullet \to X^\bullet$.  Then one defines the {\it
(strict) Greenberg} local condition at $v$ to be the induced morphism
$i_v^+ \cn C^\bullet(G_v,(X_v^+)^\bullet) \to
C^\bullet(G_v,X^\bullet)$.

For a place $v \in \Sigma'$, there is the {\it unramified} local
condition at $v$.  When $X^\bullet$ is $X$ concentrated in a degree
zero, this is $U_v^+ = C^\bullet(G_v/I_v,X^{I_v})$, with $i_v^+$ being
the inflation map.  (Our present lacking of a well-defined
Hochschild--Serre morphism makes for difficulties when $X^\bullet$ is
not concentrated in a single degree.  Since we only need the simple
case, we let the reader find a more general formulation in
\cite[Chapter 7]{nekovar}.)

The maps $i_v^+$ play the role in our situation of the subgroups
$H^1_f(G_v,X) \subseteq H^1(G_v,X)$ appearing in classical Selmer
groups.  Just as one has the notationally convenient quotients
$H^1_s(G_v,X)$, we form $U_v^- = \Cone(-i_v^+)$, which come with
canonical maps $i_v^- \cn C^\bullet(G_v,X) \to U_v^-$.  In the case of
a (strict) Greenberg local condition induced by $j_v^+ \cn X_v^+ \to
X_v$, if we put $(X_v^-)^\bullet = \Cone(-j_v^+)$, then $i_v^-$ is
induced by $j_v^-$, in the evident notation.

Now we put the local conditions together to define Selmer complexes.
We set $U_S^\pm = \bigoplus_{v \in S_f} U_v^\pm$ and $i_S^\pm =
\bigoplus i_v^\pm$.  Recall that our choices of embeddings of
algebraic closures have furnished us with maps $G_v \to G_{K,S}$.
Pulling back via these maps, we get functorial ``restriction'' maps
$\text{res}_v \cn C^\bullet(G_{K,S},X) \to C^\bullet(G_v,X)$.  We
gather these restriction maps into one map $\text{res}_{S_f} =
\bigoplus_{v \in S_f} \text{res}_v$.

Given a system $\Delta$ of local conditions, we define the {\it Selmer
complex}
\[
\wt{C}_f^\bullet(X) := \wt{C}_f^\bullet(G_{K,S},X;\Delta) :=
\Cone\left(i_S^- \circ \text{res}_{S_f}\right)[-1].
\]
We write $\wt{\bfR\Ga}_f(X)$ for its image in $\bfD(R)$, and we call
its cohomology groups $\wt{H}_f^*(X)$ the {\it extended Selmer groups}
of $X$ (in the sense of \Nekovar).

Consider our examples above of local conditions.  When we have
``imposed no conditions'' we obtain
$\wt{\bfR\Ga}_f(G_{K,S},X;\Delta_\emptyset) = \bfR\Ga(G_{K,S},X)$,
which is essentially the same as continuous \'etale cohomology over
$\Spec \calO_{K,S}$, and when we have ``imposed full conditions'' we
obtain $\bfR\Ga_c(G_{K,S},X) := \wt{\bfR\Ga}_f(G_{K,S},X;\Delta_c)$,
which is essentially continuous \'etale cohomology over $\Spec
\calO_{K,S}$ with compact supports.  For a general system $\Delta$,
one obtains from the definition of $\wt{C}_v^\bullet$ an exact
sequence
\begin{multline*}
0 \to \wt{H}_f^0(G_{K,S},X;\Delta) \to H^0(G_{K,S},X) \to \bigoplus_{v
  \in S_f} H^0(U_v^-) \to \\
\to \wt{H}_f^1(G_{K,S},X;\Delta) \to H^1(G_{K,S},X) \to \bigoplus_{v
  \in S_f} H^1(U_v^-).
\end{multline*}
This exact sequence tells us how $\wt{H}_f^1(X)$ compares to
traditional Selmer groups, and, together with global duality, how to
bound the degrees of the Selmer complex.

Given a system $\Delta$ of local conditions for $X$, we may define the
{\it $D(1)$-dual local conditions} $\Delta^*$ for $D(X^\bullet)(1)$ by
taking the exact triangles
\[
U_v^+ \to C^\bullet(G_v,X^\bullet) \to U_v^-,
\]
hitting them with $D$, twisting them by $R(1)$, applying local duality
for each $v$, and shifting the result by $[-2]$.  One similarly
obtains $\scrD(1)$-dual local conditions $\Delta^*$ for
$\scrD(T^\bullet)(1)$.

One has that $(\Delta_\emptyset)^* = \Delta_c$, and $(\Delta_c)^* =
\Delta_\emptyset$.  The $D(1)$-dual (resp.\ $\scrD(1)$-dual) to a
Greenberg local condition induced by $j_v^+ \cn X_v^+ \to X$ is the
Greenberg local condition induced by $D(j_v^-)(1)$ (resp.\
$\scrD(j_v^-)(1)$).  The $D(1)$-dual to an unramified local condition
is again the unramified local condition of the $D(1)$-dual; however,
the $\scrD(1)$-dual to an unramified local condition is {\it not
necessarily} isomorphic to the unramified local condition of the
$\scrD(1)$-dual, in general!  Assuming $T$ is concentrated in degree
$0$, one can only guarantee that unramified local conditions are
$\scrD(1)$-dual to each other when $\Frob_v-1$ acts bijectively on
$H^1(I_v,T)$.

\subsection{Global duality and Iwasawa theory}\label{sub-global}

We fix $K$ and $S$ as in the preceding section.  In this section,
$X^\bullet$ is equipped with a continuous, linear $G_{K,S}$-action,
and a choice of local conditions $\Delta$.

The global duality theorem says that $G_{K,S}$ satisfies the
finiteness hypothesis with $\text{c.d.}_p(G_{K,S}) = 2$, and that the
sum of invariants of the cup product has adjoint morphisms
\begin{eqnarray*}
\wt{\bfR\Ga}_f(G_{K,S},D(X^\bullet)(1);\Delta^*)
  & \to & D\left(\wt{\bfR\Ga}_f(G_{K,S},X^\bullet;\Delta)\right)[-3]
\qquad \text{and} \\
\wt{\bfR\Ga}_f(G_{K,S},\scrD(T^\bullet)(1);\Delta^*)
  & \to & \scrD\left(\wt{\bfR\Ga}_f(G_{K,S},T^\bullet;\Delta)\right)[-3],
\end{eqnarray*}
that are isomorphisms in $\bfD(R)$.

Consider the special case when $X$ is finite and concentrated in
degree $0$, and $\Delta = \Delta_\emptyset$, so that $\Delta^* =
\Delta_c$.  Upon projecting onto the level of cohomology, the theorem
involving $D$ says that the global Galois cohomology is Pontryagin
dual to the cohomology {\it with compact supports} of the Cartier
dual.  Thus we recover a statement reminiscent of Poincar\'e duality
from this special case.

One has a version of the above theorem that applies to limits of
cochain complexes over abelian towers of number fields.  Let $L/K$ be
a $\bbZ_p^n$-extension with Galois group $\Ga_L$, write $\La_L :=
R[\![\Ga]\!]$ as in \S\ref{sub-notation}, and let $\{K_\alpha\}$ be
the collection of subfields of $L$ that are finite over $K$.  Define
the {\it Iwasawa Selmer complexes} by the expressions
\begin{gather*}
\begin{split}
\wt{\bfR\Ga}_{f,\text{Iw}}(L/K,T;\Delta)
& := \llim_\alpha
  \wt{\bfR\Ga}_f(G_{K_\alpha,S_\alpha},T;\Delta_\alpha) \\
& \cong
  \wt{\bfR\Ga}_f(G_{K,S},T \otimes_R \La_L; \Delta \otimes_R \La_L)
\in \bfD(\La_L) \qquad \text{and}
\end{split}\\
\begin{split}
\wt{\bfR\Ga}_f(K_S/L,A;\Delta)
& := \rlim_\alpha
  \wt{\bfR\Ga}_f(G_{K_\alpha,S_\alpha},A;\Delta_\alpha) \\
& \cong
  \wt{\bfR\Ga}_f(G_{K,S},\Hom_R(\La_L,A); \Hom_R(\La_L;\Delta))
\in \bfD(\La_L),
\end{split}
\end{gather*}
where the isomorphisms are due to Shapiro's lemma, $S_\alpha$ is the
set of places of $K_\alpha$ lying over $S$, and $\Delta_\alpha$
consists of the local conditions $\Delta$, appropriately propagated to
$K_\alpha$.  (See \cite[\S\S8.6--8.7]{nekovar} for how to propagate
unramified and Greenberg conditions.)  One can show that taking
cohomology of these complexes commutes with forming the limits.  In
other words, we get
\begin{gather*}
H^*\left(\wt{\bfR\Ga}_{f,\text{Iw}}(L/K,T;\Delta)\right) \cong
\llim_\alpha \wt{H}_f^*(G_{K_\alpha,S_\alpha},T;\Delta_\alpha)
=: \wt{H}_{f,\text{Iw}}^*(L/K,T;\Delta), \\
H^*\left(\wt{\bfR\Ga}_f(K_S/L,A;\Delta)\right) \cong
\rlim_\alpha \wt{H}_f^*(G_{K_\alpha,S_\alpha},A;\Delta_\alpha)
=: \wt{H}_f^*(K_S/L,A;\Delta).
\end{gather*}

Using the above complexes, we get the Iwasawa-theoretic duality
theorem, which says that the natural maps
\begin{gather*}
\wt{\bfR\Ga}_{f,\text{Iw}}(L/K,\scrD(T^\bullet)(1);\Delta^*)
  \to \scrD_{\La_L}\left(\wt{\bfR\Ga}_{f,\text{Iw}}(L/K,T^\bullet;\Delta)\right)^\iota[-3], \\
\wt{\bfR\Ga}_f(K_S/L,D(T^\bullet)(1);\Delta^*)
  \to D_{\La_L}\left(\wt{\bfR\Ga}_{f,\text{Iw}}(L/K,T^\bullet;\Delta)\right)^\iota[-3]
\end{gather*}
are isomorphisms in $\bfD(\La_L)$.  (In light of Shapiro's lemma, this
theorem is a trivial consequence of the original global duality
theorem, but now working over the ring $\La_L$.)

\Nekovar\ has proved a general control theorem.  Let $r \in R$ be a
non-zero-divisor for all of the following: $R$, $X$, both $X_v^\pm$
(when $\Delta$ is Greenberg at $v$), and $X^{I_v}$ (when $\Delta$ is
unramified at $v$).  Moreover, when $\Delta$ is unramified at $v$,
assume that
\begin{equation}\label{eqn-hyp-control}
X^{I_v} \otimes_R R/r = (X \otimes_R R/r)^{I_v}.
\end{equation}
Then one has
\[
\wt{\bfR\Ga}_f(X) \otimes_R^\bfL R/(r) \stackrel{\sim}{\to}
\wt{\bfR\Ga}_f(X/rX).
\]
More generally, let $I \subset R$ be any ideal that is generated by a
sequence $r_1,\ldots,r_n$ that is regular for the above list of
objects and such that for $i=1,\ldots,n$ one has
\[
X^{I_v} \otimes_R R/(r_1,\ldots,r_i) = (X \otimes_R
R/(r_1,\ldots,r_i))^{I_v}.
\]
Then by induction we deduce that
\[
\wt{\bfR\Ga}_f(X) \otimes_R^\bfL R/I \stackrel{\sim}{\to}
\wt{\bfR\Ga}_f(X/IX).
\]

One knows by \cite[\S5.4, Lemma (i)]{lang} that $\bbZ_p^d$-extensions
are unramified away from $p$, including all places $v$ where $\Delta$
is the unramified condition. For every place $v'$ of $L$ lying over
such $v$, one has $I_{v'} = I_v$, and $I_v$ acts trivially on $\La_L$.
Thus $X \otimes_R \La_L$ satisfies $(X \otimes_R \La_L)^{I_v} =
X^{I_v} \otimes_R \La_L$, and, as a result, Equation
\ref{eqn-hyp-control} holds for $X$.  Thus we can apply the control
theorem to $X$ and $I = \calI = \ker(\La_L \twoheadrightarrow R)$
(with $r_i$ of the form $\bra{\ga_i}-1$ with $\ga_i \in \Ga$), to
obtain the isomorphism
\[
\wt{\bfR\Ga}_{f,\text{Iw}}(L/K,T) \otimes_\La^\bfL R
\stackrel{\sim}{\to} \wt{\bfR\Ga}_f(T).
\]
For another example, if $r \in R$ is prime then, granted Equation
\ref{eqn-hyp-control}, we can compute the Selmer and Iwasawa-theoretic
Selmer complexes of $T/rT$ by taking those associated to $T$ and
(derived-)tensoring down.  This can be used sometimes when $R$ is a
Hida--Hecke algebra and $(r) \in \Spec R$ corresponds to a classical
eigenform.

\bibliography{growth}

\end{document}